\documentclass[11pt,a4paper]{article}
\usepackage[a4paper,textwidth=147mm,textheight=218mm,centering]{geometry}
\usepackage[T1]{fontenc}
\usepackage{lmodern,amsmath,amssymb,amsthm,mathtools}
\usepackage{tikz}
\usetikzlibrary{arrows.meta,calc,positioning,patterns}
\usepackage{titlesec}
\usepackage[hidelinks]{hyperref}
\usepackage{microtype,needspace,etoolbox}
\usepackage{array}
\newcommand{\secref}[1]{Section~\ref{#1}}
\newcommand{\subsecref}[1]{Subsection~\ref{#1}}
\newcommand{\secrefs}[2]{Sections~\ref{#1} and~\ref{#2}}
\titleformat{\subsection}[runin]{\normalfont\bfseries}{\thesubsection}{.5em}{}[.]
\titlespacing*{\subsection}{0pt}{1.6ex plus .4ex minus .2ex}{.6em}
\newtheorem{theorem}{Theorem}[section]
\newtheorem{maintheorem}{Theorem}

\newtheorem{proposition}[theorem]{Proposition}
\newtheorem{lemma}[theorem]{Lemma}
\newtheorem{corollary}[theorem]{Corollary}
\theoremstyle{definition}

\theoremstyle{remark}

\newcommand{\src}{\mathsf{s}}

\newcommand{\C}{\mathbb Z/2\mathbb Z}
\newcommand{\norm}[1]{\lVert#1\rVert_{\mathbb R/\mathbb Z}}
\definecolor{conn}{RGB}{0,92,89}
\definecolor{bound}{RGB}{151,82,25}
\tikzset{edge/.style={-{Latex[length=2mm]},line width=.7pt},
  connector/.style={conn,<->,>=Latex,line width=.8pt},
  point/.style={circle,fill=black,inner sep=1.2pt},
  pending/.style={circle,draw=bound,fill=white,line width=.7pt,inner sep=1.5pt},
  every picture/.style={font=\small}}
\title{Groups of intermediate growth arising from irrational rotations}
\author{Zheng Kuang\thanks{School of Mathematics, South China University of Technology, Guangzhou, China; Email: \href{mailto:kzkzkzz@scut.edu.cn}{kzkzkzz@scut.edu.cn}; ORCID: \href{https://orcid.org/0000-0001-6581-8207}{0000-0001-6581-8207}.}}
\date{}
\begin{document}
\maketitle
\begin{abstract}
For each irrational angle $\theta$, we fragment the associated minimal dihedral group action on the split circle to obtain a periodic group generated by five involutions. We give a sufficient condition for intermediate growth of the group that holds, in particular, when the denominators $q_n(\theta)$ of the convergents of $\theta$ satisfy $\log q_n(\theta)=O(n^d)$ for some $d>0$. If $\log q_n(\theta)=O(n)$, then the growth function of the group is at most $\exp(CR^\alpha)$ for some $\alpha<1$. We bound the numbers of traverses and portraits by hitting estimates for intervals of length comparable to $1/q_n$. If the partial quotients are unbounded, then the orbital graphs at points with trivial stabilizer in the dihedral group action are not linearly repetitive. For almost every angle, the group has intermediate growth with this bound and these graphs are not linearly repetitive. We realize the fragmented action by a tile inflation process over the Bratteli diagram described by the author and express the conditions using return times on the modular surface.
\end{abstract}
\noindent\textbf{Keywords.} groups of intermediate growth; irrational rotations; Schreier graphs; continued fractions; tile inflation; AF-by-discrete groupoids

\noindent\textbf{2020 Mathematics Subject Classification.} Primary 20F69, 37E10; Secondary 22A22, 37B10, 37D40, 11J70.
\enlargethispage{\baselineskip}
\vspace{-\baselineskip}
\begingroup
\makeatletter
\patchcmd{\l@section}{\addvspace{1.0em \@plus\p@}}{\addvspace{.65em \@plus\p@}}{}{}
\makeatother
\tableofcontents
\endgroup

\section{Introduction}
We study groups obtained by fragmenting the non-free minimal actions of the infinite dihedral group $D_\infty$ associated with irrational rotations, using the Bratteli diagrams and tile inflation processes described in \cite[Section~5]{Kua25}. We extend the method of traverses of \cite{BNZ22} to obtain subexponential growth estimates when the orbital graphs at points with trivial stabilizer in the dihedral group action are not linearly repetitive. We also relate these actions to returns of the geodesic flow on the modular surface through first return maps of the translation flows for lattice strips.

Milnor's question~\cite{Mil68} as to whether the growth of every finitely generated group is either polynomial or exponential was answered in the negative by Grigorchuk \cite{Gri84}. Nekrashevych \cite{Nek18} developed a construction of periodic groups by fragmenting minimal dihedral group actions. In this construction, one replaces a reflection by several homeomorphisms, each of which agrees pointwise with the reflection or the identity. The Grigorchuk groups \cite{Gri80,Gri84} arise from fragmentations of the dihedral group action given by the binary Gray code (see \cite{Vor12,Nek18}).

The theorems on growth in \cite[Theorem~3.1]{BNZ22} and \cite[Theorem~2.9]{Nek25} (the latter for admissible subshifts) assume linear repetitivity. For the linear orbital graphs considered here, linear repetitivity bounds the gaps between occurrences of each labeled segment by a constant times its length, not merely the number of distinct segments of that length. We instead estimate the occurrences of particular labeled segments.

Let $\theta\in(0,1)$ be irrational. We work on the split circle $\mathcal{X}_\theta$, on which the reflections $a(x)=-x$ and $b(x)=\theta-x$ lift to homeomorphisms and $ba$ acts as a rotation. The dihedral group $D_{\infty}$ is generated by $\{a,b\}$. The notions used below are recalled in \secref{sec:prelim}.

In \cite[Theorem~5.2]{Kua25}, we conjugated the action of $D_\infty$ on $\mathcal{X}_\theta$ to an action on the space of infinite paths $\Omega(\mathsf{B})$ of a Bratteli diagram $\mathsf B$ constructed from the continued fraction of $\theta$. See \cite[Subsection~5.2]{Kua25} for the tile inflation process defining the action $D_\infty\curvearrowright\Omega(\mathsf{B})$ and \cite[Subsection~2.1]{Kua26} for the general construction. In this realization, the action is linearly repetitive and of bounded type if and only if the partial quotients are bounded (\cite[Proposition~5.5 and Corollary~5.6]{Kua25}). The growth estimate for the fragmentation in \cite[Subsection~5.3]{Kua25} requires, in addition, partial quotients greater than one and uniformly bounded gaps between odd partial quotients. These restrictions control when certain edges labeled by $a$ are introduced in the successive subdivisions (see \cite[Proposition~5.7]{Kua25}).

We retain these diagrams but fragment only $a$, whereas in~\cite{Kua25}, both $a$ and $b$ were fragmented. In the orbital graphs of the dihedral group action, moving along an edge labeled by $a$ and then one labeled by $b$ amounts to applying the rotation $ba$. We choose the clopen sets defining the fragments using hitting and separation estimates for this rotation. Thus, we control the occurrences of the prescribed labeled segments used to count traverses and elements of portraits for arbitrary words in the generators of the group obtained by fragmentation. Our hitting estimates for shrinking intervals are valid for every initial point, in contrast to the metric results on shrinking targets for rotations in \cite{Kur55,Kim07,Tse08,CC18}. Instead of bounding the partial quotients, we impose a condition on the denominators $q_n(\theta)$. We prove the following (Theorem~\ref{thm:growth}).

\Needspace{6\baselineskip}
\begin{maintheorem}\label{main:growth}
For each irrational $\theta\in(0,1)$, the fragmentation of the associated minimal dihedral group action constructed in \secref{sec:construction} gives an infinite periodic group $G_\theta$ generated by $5$ involutions. Let $\gamma_{G_\theta}$ be its growth function with respect to these generators, and let $q_n(\theta)$ be the denominators of the continued-fraction convergents of $\theta$. If
\[
 \liminf_{n\to\infty}(2/3)^{n/64}\log q_n(\theta)=0,
\]
then $G_\theta$ has intermediate growth. Put
\[
 \Lambda_\theta=\limsup_{n\to\infty}\frac{\log q_n(\theta)}{n}.
\]
If $\Lambda_\theta<\infty$, then for every
\[
 \frac{64\Lambda_\theta}{64\Lambda_\theta+\log(3/2)}<\alpha<1,
\]
there exists $C>0$ such that $\gamma_{G_\theta}(R)\leq\exp(CR^\alpha)$. For every irrational $\theta$ and every point with trivial stabilizer in the dihedral group action, the number of isomorphism classes of rooted labeled balls of radius $r$ in the orbital graph of the fragmented action at this point is $O(r+1)$. If the partial quotients of $\theta$ are unbounded, then the orbital graphs at points with trivial stabilizer in the dihedral group action are not linearly repetitive.
\end{maintheorem}

By L\'evy's theorem \cite{Lev36}, $\log q_n(\theta)/n\to\pi^2/(12\log 2)$ for almost every $\theta$. The partial quotients are unbounded almost everywhere by the Borel--Bernstein theorem (see \cite{Khi64}). Thus, for almost every angle, the group in Theorem~\ref{main:growth} has intermediate growth, while the orbital graphs of the fragmented action at points with trivial stabilizer in the dihedral group action are not linearly repetitive. The assertion on repetitivity concerns these particular actions. In \secref{sec:repetitivity}, we give an example with unbounded partial quotients, only finitely many odd partial quotients, and $\Lambda_\theta=\log(1+\sqrt2)$.

See also \cite{MB15} for embeddings into topological full groups of minimal subshifts. Subexponential growth of such groups need not give a bound of the form $\exp(CR^\alpha)$ with $\alpha<1$, as shown by Erschler \cite{Ers04} for certain Grigorchuk groups. For $\theta=(\sqrt{5}-1)/2$, the reflections of the circle are those used in the construction of Nekrashevych in \cite[Section~8]{Nek18}, with $a$ and $b$ interchanged. Bartholdi, Nekrashevych and Zheng \cite[Subsection~4.2]{BNZ22} obtained the estimate $\gamma(R)\preccurlyeq\exp(R^{0.9181})$ for the group constructed there. We use a different fragmentation and do not seek optimal exponents in Theorem~\ref{main:growth}.

\subsection{Sketch of the method}
The method of \cite{BNZ22} relates traverses to standard portraits, using linear repetitivity to find prescribed labeled segments. We choose the fragments of $a$ so that the occurrences follow from the rotation itself.

Choose $Q_j=q_{n_j}(\theta)$ with $n_{j+1}-n_j=32$. The fragmentation is defined by three pairs of intervals near the fixed point $\eta=1/2$ of $a$, of lengths comparable to $1/Q_j$. Every block of $Q_j$ consecutive rotation iterates meets each prescribed interval on one side of $\eta$. Separation estimates identify the two branches at the corresponding $a$-edges with the required labeled segments of the orbital graph at $\eta$. Neither estimate requires restrictions on the sizes or parities of the partial quotients. The counting argument of \cite[proof of Theorem~3.1]{BNZ22} gives $T_{j+1}(w)\leq(2/3)T_j(w)$ for the numbers $T_j(w)$ of return paths counted in \secref{sec:contraction}.

In \cite[Lemma~2.4]{BNZ22}, linear repetitivity bounds the number of isomorphism classes of rooted labeled balls. The proof of \cite[Proposition~2.7]{BNZ22} uses the gap condition of \cite[Subsection~2.2]{BNZ22} and bounded ratios of the lengths of successive segments to bound the multiplicity of the correspondence between portraits and traverses. Instead, in Lemma~\ref{lem:complexity}, we bound the number of isomorphism classes of labeled neighborhoods of radius $2R+2$ of the selected $a$-edges by $O(1+R/Q_j)$ for $R\geq Q_j$. Let $\mathcal P_w$ be the standard portrait of a word $w$ of length at most $R$ (see \cite[Definition~2.2]{BNZ22}). Counting separately trajectories of radius less than $16Q_j$ and trajectories crossing one of the prescribed segments gives, in Lemma~\ref{lem:portrait},
\[
 |\mathcal P_w|\leq C\left(Q_j+\frac{R^2}{Q_j}(2/3)^j\right),\qquad R\geq16Q_j,\quad j\geq2.
\]
Taking $R_j=\lceil Q_j(3/2)^{j/2}\rceil$ gives $\log\gamma_{G_\theta}(R_j)\leq CQ_j\log R_j$. Dividing by $R_j$ and using submultiplicativity gives subexponential growth under the first condition of Theorem~\ref{main:growth}. If $\log Q_j=O(j)$, then we also obtain the stronger upper bound. No uniform bound on the ratios of the lengths of successive segments is needed.

\subsection{Tile inflation and the modular surface}
The tile inflation processes of \cite{Kua25,Kua26a} describe actions on the spaces of infinite paths of Bratteli diagrams. The tail groupoid records the movements within finite tiles. In our realization, the remaining germs are obtained from germs at fixed points of the reflections. We realize the fragmentation over a telescoping of the diagram of \cite[Subsection~5.1]{Kua25}, using the conventions of \subsecref{subsec:prelim-tiles}. The following theorem is proved in \secref{sec:native} (Theorem~\ref{thm:hierarchy}).

\begin{maintheorem}\label{main:hierarchy}
For every irrational $\theta\in(0,1)$, the action of $G_\theta$ on the split circle $\mathcal{X}_{\theta}$ is conjugate to an action defined by a tile inflation process over a telescoping of the Bratteli diagram $\mathsf B$ constructed in \cite[Theorem~5.2]{Kua25} for the associated dihedral group action. The tail groupoid $\mathfrak T_{\mathsf B}$ in the realization of the dihedral group action is an open subgroupoid of the groupoid of germs of $G_\theta$, with discrete complement.
\end{maintheorem}
If the partial quotients are unbounded, then the numbers of edges at the levels of every telescoping of $\mathsf B$ are unbounded (see \secref{sec:native}).

Series \cite{Ser85} related cutting sequences of geodesics on the modular surface to continued fractions, and the coding of the geodesic flow on two sheets is due to Arnoux \cite{Arn94} (see also \cite{AS13,AS26}). The scenery flow of Arnoux and Fisher \cite{AF01,AF05} describes subdivision and rescaling of tilings associated with Sturmian systems. In \secrefs{sec:strip}{sec:returns}, we associate to a lattice strip the action groupoid of the first return map of the translation flow to the canonical transversal, with its displacement cocycle and a distinguished unit. Proposition~\ref{main:reconstruction} recovers the lattice and width of the strip from these data. For a reduced horizontal strip with rotation parameter $y$, this first return map is conjugate to a rotation of $\mathcal{X}_y$. Its successive inductions describe modular returns by restrictions to clopen orbit transversals, with the induced cocycles and images of the distinguished units; we retain the orientation sheet and suspension coordinate. Adjoining the reflection $a$ gives the dihedral group action of \cite[Theorem~5.2]{Kua25} for $\theta=y$, from which we construct $G_\theta$.

The partial quotients of $y$ determine both the subdivisions of tiles associated with backward modular returns and the incidence matrices of the Bratteli diagram used for the dihedral group action. By Lemma~\ref{lem:boundary}, the endpoints of edges introduced by connectors labeled by $a$ at higher levels lie in shrinking neighborhoods of the persistent boundary point $\eta$ (a boundary point of an infinite tile). Thus, we choose the intervals defining the fragments of $a$ near $\eta$ and show in \secref{sec:native} that the corresponding edges are introduced at the prescribed levels of a telescoping. The denominators $q_n(y)$ of the convergents of the same continued fraction control the lengths of these intervals and the gaps between occurrences of the edges in the orbital graphs. The separation estimates determine the labels in the prescribed neighborhoods of these edges. In \secrefs{sec:contraction}{sec:portraits}, we use the resulting labeled segments to obtain the contraction estimate for the numbers of return paths and the bound on the cardinalities of portraits, respectively, and hence bounds on the growth function of the group.

The successive backward modular returns determine clopen orbit transversals $C_N^-$ with
\[
 -\log\mu_y(C_N^-)=\log q_N(y)+O(1)
 =\frac12\sum_{n=-N}^{-1}\tau_n+O(1),
\]
where $\mu_y$ is the invariant probability measure and $\tau_n$ are the Gauss return times. Exchanging the coordinate directions with a sign change gives the analogous comparison for forward returns. This relates the conditions in Theorem~\ref{main:growth} to modular returns (Corollary~\ref{cor:modular}).

\subsection{Organization of the paper}
After recalling the necessary notions in \secref{sec:prelim}, we describe the reconstruction of lattices and modular returns in \secrefs{sec:strip}{sec:returns}. We construct the fragmentation in \secref{sec:construction} and prove Theorem~\ref{main:hierarchy} in \secref{sec:native}. In \secrefs{sec:contraction}{sec:portraits}, we establish the estimates on traverses and portraits, respectively. We complete the proof of Theorem~\ref{main:growth} in \secref{sec:repetitivity}.

The proof of Theorem~\ref{main:growth} does not require the results in \secrefs{sec:strip}{sec:returns} or the realization in \secref{sec:native}. A reader interested primarily in growth may therefore read Subsections~\ref{subsec:prelim-growth}--\ref{subsec:prelim-traverses} and the first paragraph of \subsecref{subsec:prelim-modular}, followed by \secref{sec:construction} and Sections~\ref{sec:contraction}--\ref{sec:repetitivity}. A reader interested in the modular surface may read \secrefs{sec:strip}{sec:returns} using the statement of Theorem~\ref{main:growth}; the construction of the group to which the modular conditions apply is given in \secref{sec:construction}.

\section{Preliminaries}\label{sec:prelim}
\subsection{Growth of groups and orbital graphs}\label{subsec:prelim-growth}
Let $G$ be a group generated by a finite symmetric set $\mathcal S$. We use left actions and compose maps from right to left. The length $|g|$ of an element $g\in G$ is the least length of a word in $\mathcal S$ representing $g$. The growth function is
\[
 \gamma_G(R)=\#\{g\in G:|g|\leq R\}.
\]
We say that $G$ has \emph{subexponential growth} if $R^{-1}\log\gamma_G(R)\to0$. If, in addition, its growth function is not bounded above by any polynomial, then $G$ has \emph{intermediate growth}. A group is said to be \emph{periodic} if each of its elements has finite order. By Gromov's theorem \cite{Gro81}, a finitely generated group of polynomial growth is virtually nilpotent. If it is also periodic, then it is finite. Thus, every infinite finitely generated periodic group of subexponential growth has intermediate growth.

Suppose that $G$ acts on a set $\mathcal{X}$. The orbital graph at $\zeta\in\mathcal{X}$ is the graph with vertices $G\zeta$ and an edge from $v$ to $s(v)$ labeled by $s$ for each $v\in G\zeta$ and $s\in\mathcal S$. We retain loops and edges with different labels and measure distances in the underlying undirected graph. An orbital graph is said to be \emph{linear} if it is quasi-isometric to the real line or a half-line. The trajectory of a word $w=s_m\cdots s_1$ starting at $v$ is the sequence of vertices obtained by applying $s_1,\ldots,s_m$ in that order. If a group of homeomorphisms has a dense orbit, then its action on that orbit is faithful.

A labeled orbital graph is said to be \emph{linearly repetitive} if there is $L\geq1$ such that, for every $r\geq1$, every ball of radius $\lceil Lr\rceil$ contains a copy of every rooted labeled ball of radius $r$ in the graph. If the underlying simple graph is a line, then this is equivalent to the existence of $L'\geq1$ such that, for every integer $r\geq1$, every segment of length $\lceil L'r\rceil$ contains a copy of every segment of length $r$. A copy is given by a label-preserving isomorphism, which may reverse the orientation of the line. See \cite[Subsection~2.2]{BNZ22}.

\subsection{Irrational rotations and fragmentations}\label{subsec:prelim-frag}
For an irrational $\theta\in(0,1)$, write $\mathcal{X}_\theta$ for the split circle obtained from $\mathbb R/\mathbb Z$ by replacing every point of $\mathbb Z\theta\pmod1$ by two one-sided lifts (see the beginning of~\cite[Subsection 5.1]{Kua25} for more detail). An arc with endpoints in this orbit gives a clopen interval by taking the right lift of its left endpoint and the left lift of its right endpoint. These intervals form a basis of the topology of the Cantor set $\mathcal{X}_\theta$. We denote by the same symbols the two reflections on the circle and their lifts:
\[
 a(x)=-x,\qquad b(x)=\theta-x\pmod1.
\]
Since the reflections reverse orientation, their lifts interchange the two sides at split points. Their product $ba$ is the lift of the rotation by $+\theta$. This gives a minimal action of $D_\infty=\langle a,b\mid a^2=b^2=1\rangle$. The orbital graphs of this action are linear. At a point with trivial stabilizer, the orbital graph is a line, with edges labeled alternately by $a$ and $b$. At a fixed point of either reflection, it is a ray with a loop at its endpoint.

The rotation of $\mathcal{X}_\theta$ has a unique invariant probability measure $\mu_\theta$, since every invariant probability measure projects to Lebesgue measure and the countable set of split points has measure zero. If a Borel subset of $\mathcal{X}_\theta$ has boundary of measure zero, then the frequency of visits along every orbit of the rotation equals its measure.

Following \cite{Nek18}, a \emph{fragmentation} of an involution $s$ of a Cantor set is a finite group of homeomorphisms $A$ such that $u(\zeta)\in\{\zeta,s(\zeta)\}$ for each $u\in A$, and, for every $\zeta$, there is $u\in A$ with $u(\zeta)=s(\zeta)$. Since $A$ and $\langle s\rangle$ have the same orbits, fragmenting the generators of a dihedral group action preserves the orbits but changes the labels of the edges. In this paper, we fragment $a$ and retain $b$. Let $A$ and $B$ be fragmentations of the generators $a$ and $b$, respectively, of a minimal dihedral group action. If $a(\eta)=\eta$ and the interior of the set of fixed points of each $u\in A$ accumulates at $\eta$, then $\langle A\cup B\rangle$ is periodic (\cite[Theorem~1.1]{Nek18}). The group is infinite, since its orbits are those of the minimal dihedral group action on the Cantor set. We use this criterion in Proposition~\ref{prop:periodic}.

\subsection{Traverses and portraits}\label{subsec:prelim-traverses}
We recall the conventions of \cite[Subsections~2.1--2.2]{BNZ22}, using left actions. In a linear orbital graph whose underlying simple graph is a line or a ray, a segment is the subgraph induced by a finite interval of vertices. Its length is the number of edges between consecutive distinct vertices. Choose a direction on a nontrivial segment, so that its initial and terminal vertices are, respectively, the entrance and the exit. A \emph{traverse} of this segment is a portion of the trajectory of a word that goes from the entrance to the exit, stays in the segment and meets neither endpoint at any intermediate time. We also call such a traverse a \emph{complete crossing}. If a trajectory visits the entrance and subsequently the exit, then take the first visit to the exit after the specified visit to the entrance. The portion of the trajectory ending at this arrival and starting at the last preceding visit to the entrance is a traverse. Occurrences of a subword at different starting times are counted separately. In \secref{sec:contraction}, we identify two branches of a segment with a segment of the orbital graph at a fixed point of a reflection. The image of a traverse then leaves the outer endpoint and returns to it; we count these return paths, subject to an additional condition on their steps at the endpoint of the ray.

Let $w=s_m\cdots s_1$ and let $B_v(r)$ be the ball of radius $r$ about a vertex $v$ of an orbital graph. Denote by $r_w(v)$ the maximum distance from $v$ attained by the trajectory of $w$. The \emph{standard portrait} $\mathcal P_w$ of $w$ (\cite[Definition~2.2]{BNZ22}) is the set of isomorphism classes of the bi-rooted labeled graphs
\[
 (B_v(r_w(v)),v,w(v)),\qquad v\in G\zeta.
\]
The second root marks the final vertex. For a linear orbital graph whose underlying simple graph is a line, we retain the orientations of the balls. At every vertex $v$, at least one graph in the portrait admits an embedding sending its first root to $v$. The labels determine the trajectory of $w$, so every such embedding sends the second root to $w(v)$. Thus, if two words have the same portrait on an orbit on which the action is faithful, then they represent the same group element. To bound the growth function, it suffices to bound the number of possible portraits of words of length at most $R$. This is the counting method of \cite[Proposition~2.3]{BNZ22}.

\subsection{Groupoids}\label{subsec:prelim-groupoids}
A groupoid $\mathfrak{G}$ is a category in which every arrow has an inverse. For an action of a discrete group $G$ on a Cantor set $\mathcal{X}$, the action groupoid $\mathcal{X}\rtimes G$ has arrows $(g,\zeta)$ from $\zeta$ to $g(\zeta)$ and composition $(h,g(\zeta))(g,\zeta)=(hg,\zeta)$. We write $\mathsf{s}$ and $\mathsf{r}$ for the source and range maps. In the groupoid of germs, two arrows with the same source are identified if and only if the corresponding maps agree on a neighborhood of that source. We denote a germ by $(g,\zeta)$ as well. The arrows from a unit to itself form its isotropy group. At a point fixed by a group element, its germ can be nontrivial even though the corresponding edge in the orbital graph is a loop. An action is \emph{topologically free} if the set of fixed points of every nonidentity element has empty interior. In this case, its action groupoid and groupoid of germs are isomorphic.

A \emph{bisection} is a subset of $\mathfrak{G}$ on which both source and range maps are injective. An \emph{\'etale groupoid} is a topological groupoid for which these maps are local homeomorphisms. The topologies of the action groupoid and the groupoid of germs have bases of open bisections obtained by restricting group elements to clopen subsets. A groupoid is said to be \emph{minimal} if all its orbits are dense. If $C$ is a clopen subset of the unit space of $\mathfrak G$, then $\mathfrak G|_C$ consists of the arrows with source and range in $C$. We call $C$ an \emph{orbit transversal} if it meets every orbit. The restriction to such a clopen set is Morita equivalent to $\mathfrak G$. For a minimal homeomorphism of a Cantor set, the first-return transformation to any nonempty clopen subset generates the restriction of its action groupoid.

A \emph{real-valued cocycle} is a continuous function $c$ on a groupoid with $c(gh)=c(g)+c(h)$ whenever the composition is defined. We use such a cocycle to record the signed displacements between projected points of a lattice strip. A distinguished unit specifies the pointed configuration from which these displacements are measured.

A groupoid is said to be \emph{principal} if its isotropy groups are trivial. An elementary subgroupoid is a compact open principal subgroupoid with finite orbits. An \emph{AF groupoid} on a Cantor set is an increasing union of elementary subgroupoids, each with the same unit space. An \emph{AF-by-discrete groupoid} on a Cantor set contains an open AF subgroupoid with the same unit space and with discrete complement (see \cite[Subsection~5.2.4]{Nek22}). In our realization, movements within finite tiles give the AF subgroupoid. The groupoid of germs retains information not determined by the orbital graph.

\subsection{Bratteli diagrams and tile inflation}\label{subsec:prelim-tiles}
We recall the conventions of \cite[Section~2]{Kua25} in the case needed in this paper. A \emph{Bratteli diagram} $\mathsf B$ is a graded directed graph with finite sets of vertices $V_n$ on level $n$ and finite sets of edges $E_n$ from level $n-1$ to level $n$. There is one vertex at level zero. The space $\Omega(\mathsf B)$ of infinite paths has a basis consisting of clopen cylinders defined by finite prefixes. Two infinite paths are related by the tail groupoid $\mathfrak T_{\mathsf B}$ if and only if they agree from some level on. Its topology is given by bisections replacing a finite prefix by another ending at the same vertex. At each level, these replacements form an elementary subgroupoid, and their union is AF. Telescoping replaces paths between successive selected levels by single edges, preserving the path space and the tail groupoid.

In a \emph{tile inflation process}, the vertices of the finite labeled tile associated with $v\in V_n$ are the paths from level zero to $v$. Each edge $e$ of the diagram from $v$ to $v'\in V_{n+1}$ defines an embedding of this tile into the tile associated with $v'$, sending $p$ to $pe$. A \emph{graph with boundary} has, in addition to its internal edges, edges with one endpoint not in the set of vertices. These are called boundary edges, and the vertices to which they are attached are called boundary vertices. Two boundary edges are said to be \emph{compatible} if they have the same label and one has no source while the other has no range. The set of connectors on each level specifies which pairs of compatible boundary edges of the copies of tiles on the previous level are to be joined. In the construction below, we also complete boundary edges to give loops at vertices fixed by the corresponding generators. The infinite tile of $x\in\Omega(\mathsf B)$ is the inductive limit of the finite tiles containing its prefixes under these embeddings. An infinite path is a \emph{boundary point} of its infinite tile (or \emph{persistent boundary point}) if all its finite prefixes are boundary vertices of the corresponding finite tiles. Write $x=p\omega$, where $p$ is a finite prefix of $x$. If there is an internal edge from $p$ to $q$ labeled by $s$ in the corresponding finite tile, then we define $s(x)=q\omega$. Since the embeddings preserve the internal edges and their labels, this definition does not depend on the level. We require that the resulting partial homeomorphisms extend continuously to homeomorphisms of $\Omega(\mathsf B)$, including at the boundary points. In \secref{sec:native}, we label the edges of the finite path graphs of the original action to realize the fragmented action. A tile inflation process is of \emph{bounded type} if the numbers of boundary vertices of its finite tiles are uniformly bounded, the entire path space contains only finitely many boundary points of infinite tiles, and the numbers of vertices and edges at each level of the diagram are uniformly bounded (\cite[Definition~2.13]{Kua25}). We do not assume a uniform bound on the number of edges at the levels in \secref{sec:native}.

In our construction, the Bratteli diagram is obtained from nested Rokhlin--Kakutani partitions of the split circle. A tower consists of disjoint clopen sets, called \emph{floors}, and specified maps between consecutive floors. One of the floors is the base of the tower. The floors refine those at the previous level, and each tower is assembled from copies of lower-level towers. If the partitions generate the topology, then nested floors identify the Cantor set with the path space. The incidence matrices count the copies of each tile in the tiles at the next level.

\subsection{Continued fractions and modular returns}\label{subsec:prelim-modular}
For an irrational number $x\in(0,1)$, the Gauss map $T(x)=1/x-\lfloor1/x\rfloor$ removes the first partial quotient of its continued fraction. We write $\theta=[0;c_1,c_2,\ldots]$ and $p_n/q_n$ for the convergents of $\theta$, with $q_0=1$ and $q_{-1}=0$. The denominators satisfy $q_n=c_nq_{n-1}+q_{n-2}$. The sizes of $q_n$ govern the distribution of the first $q_n$ iterates of a rotation, as recalled in Lemma~\ref{lem:arithmetic}.

The modular surface is the quotient of the upper half-plane $\mathbb H$ by $\mathrm{PSL}_2(\mathbb Z)$ acting by fractional linear transformations. Its unit tangent bundle can be described, with an appropriate choice of coordinates, as the space $\mathrm{SL}_2(\mathbb R)/\mathrm{SL}_2(\mathbb Z)$ of unimodular lattices (lattices in $\mathbb R^2$ of covolume one). In the coordinates used here, geodesic flow is represented by the diagonal deformation $A_t=\operatorname{diag}(e^{t/2},e^{-t/2})$. A cross-section is a set of points at which we record the orbits of the flow. The first-return map sends each point to the next intersection of its forward orbit with the cross-section. The elapsed time is the return time, which defines the roof function on the cross-section. The suspension of the return map with this function recovers the flow on the part under consideration. The Farey tessellation has geodesic edges between rational endpoints $p/q$ and $r/s$ (including infinity) whenever $|ps-qr|=1$. Continued fractions describe the cutting sequences of geodesics with respect to this tessellation (\cite{Ser85}).

The Gauss map is not invertible because it discards a partial quotient. Its natural extension retains the past in a second coordinate $y$: the continued fraction of $x$ contains the future digits, and that of $y$ contains the past digits in reverse order. The continued fraction of the new past coordinate begins with $a_0=\lfloor1/x\rfloor$, followed by the digits of $y$. In the cross-section for the geodesic flow, we also retain the orientation sheet, which changes at each return (\cite{Arn94,AS13,AS26}). The precise formulas and the two conventions for return times are given in \subsecref{subsec:induced}. The description of the geodesic flow on the space of unimodular lattices allows us to describe these returns by inducing rotations on split circles associated with lattice strips.

A shrinking-target question for a rotation concerns visits of its iterates to intervals whose lengths tend to zero. In \subsecref{subsec:ranks}, we bound hitting times for intervals of length comparable to $1/q_n$ for every initial point. The intervals used to construct the fragments are distinct from those obtained by successive restrictions in \secref{sec:returns}.

\section{Groupoids associated with lattice strips}\label{sec:strip}
Let $p_1,p_2:\mathbb R^2\to\mathbb R$ be the coordinate projections, and let $\Lambda$ be a unimodular lattice such that both projections are injective on $\Lambda$ and $p_2(\Lambda)$ is dense. For $w>0$, put
\[
 P(\Lambda,w)=p_1\bigl(\Lambda\cap(\mathbb R\times[0,w))\bigr).
\]
Let $\mathfrak G_{\Lambda,w}$ be the action groupoid of the first return map of the translation flow to its canonical transversal, the closure of $\{P(\Lambda,w)-p:p\in P(\Lambda,w)\}$, with signed displacement cocycle $c$ and distinguished unit $\xi=P(\Lambda,w)$ rooted at $0$. For reconstruction of cut-and-project schemes and windows from model-set dynamical systems, see \cite[Sections~4--5]{BLM07}.

\begin{proposition}\label{main:reconstruction}
The groupoid $\mathfrak G_{\Lambda,w}$ is minimal and admits a presentation as an AF-by-discrete groupoid. The triple $(\mathfrak G_{\Lambda,w},c,\xi)$ determines $(\Lambda,w)$ in the specified physical and internal coordinates. More precisely, set $P=\{c(g):\src(g)=\xi\}$, $D=\langle P\rangle_{\mathbb Z}$ and $D^+=\langle P\rangle_{\mathbb N\cup\{0\}}$. Then $D=p_1(\Lambda)$ and the order on $D$ defined by $D^+$ is the order induced by the internal projection. Unimodularity determines the internal projection from this order, and the width of the strip is recovered from $P$.
\end{proposition}
Here, isomorphisms of triples preserve the distinguished unit and the exact real-valued cocycle.

\subsection{The first return map of the translation flow}
Put $P=P(\Lambda,w)$. For each nonempty open interval $J\subset[0,w)$, the set $p_1(\Lambda\cap(\mathbb R\times J))$ is relatively dense in $\mathbb R$. Fix a compact fundamental parallelogram for $\Lambda$. Every point $(t,0)$ differs from a lattice point by a vector in this parallelogram, and thus the internal coordinates of these lattice points lie in a bounded interval. Since $p_2(\Lambda)$ is dense, finitely many translates of $J$ by internal coordinates of lattice vectors cover this bounded interval. Adding one of the corresponding lattice vectors gives a point of the strip over $J$ whose physical coordinate is at uniformly bounded distance from $t$.

Since each compact rectangle contains finitely many lattice vectors and both projections are injective, $P$ is uniformly discrete. By relative density, enumerate $P$ as $(p_i)_{i\in\mathbb Z}$ in increasing order, with $p_0=0$. The positive differences $p_{i+1}-p_i$ belong to a finite set, since they come from lattice vectors in a bounded rectangle. For the same reason, there are only finitely many possible configurations $(P-p_i)\cap[-R,R]$ for any $R>0$.

The translation hull of $P$ is the closure of its real translates in the topology of agreement on large bounded intervals up to small translations. The canonical transversal $\mathcal{X}_{\Lambda,w}$ of the flow $\xi\mapsto\xi-t$ consists of the configurations containing $0$, i.e., the closure of $\{P-p_i:i\in\mathbb Z\}$. Equivalently, $\mathcal{X}_{\Lambda,w}$ is identified with the subshift given by the orbit closure of the sequence $(p_{i+1}-p_i)_{i\in\mathbb Z}$. The first return map $S$ of the flow to $\mathcal{X}_{\Lambda,w}$ moves the root to the next point on the right and is a homeomorphism of $\mathcal{X}_{\Lambda,w}$.

The internal coordinates of the points of $P$ at which a given finite configuration occurs are determined by whether they belong to finitely many translates of $[0,w)$. At an endpoint, the conditions for membership agree with their limits from the right. Hence, each configuration is realized over an interval of positive length and occurs with bounded gaps by relative density. Thus, $S$ is minimal. If $P+d=P$ with $d\ne0$, then $d\in p_1(\Lambda)$, and translation by $d$ adds a constant internal coordinate to the set of internal coordinates of the points in $P$. These coordinates are dense in $[0,w)$, so this constant must be zero, contradicting injectivity of $p_2$. It follows that the subshift is infinite and the action is free. Thus, $\mathcal{X}_{\Lambda,w}$ is a Cantor set.

For $\xi\in \mathcal{X}_{\Lambda,w}$, let $c((S,\xi))$ be the first return time of the translation flow, i.e., the positive distance from the root of $\xi$ to the next point on the right. Extend this function by the cocycle identity to the action groupoid $\mathfrak G_{\Lambda,w}=\mathcal{X}_{\Lambda,w}\rtimes_S\mathbb Z$; in particular, the displacements for negative powers of $S$ are negative. We use $(S^i,\xi)$ for the arrow from $\xi$ to $S^i(\xi)$.

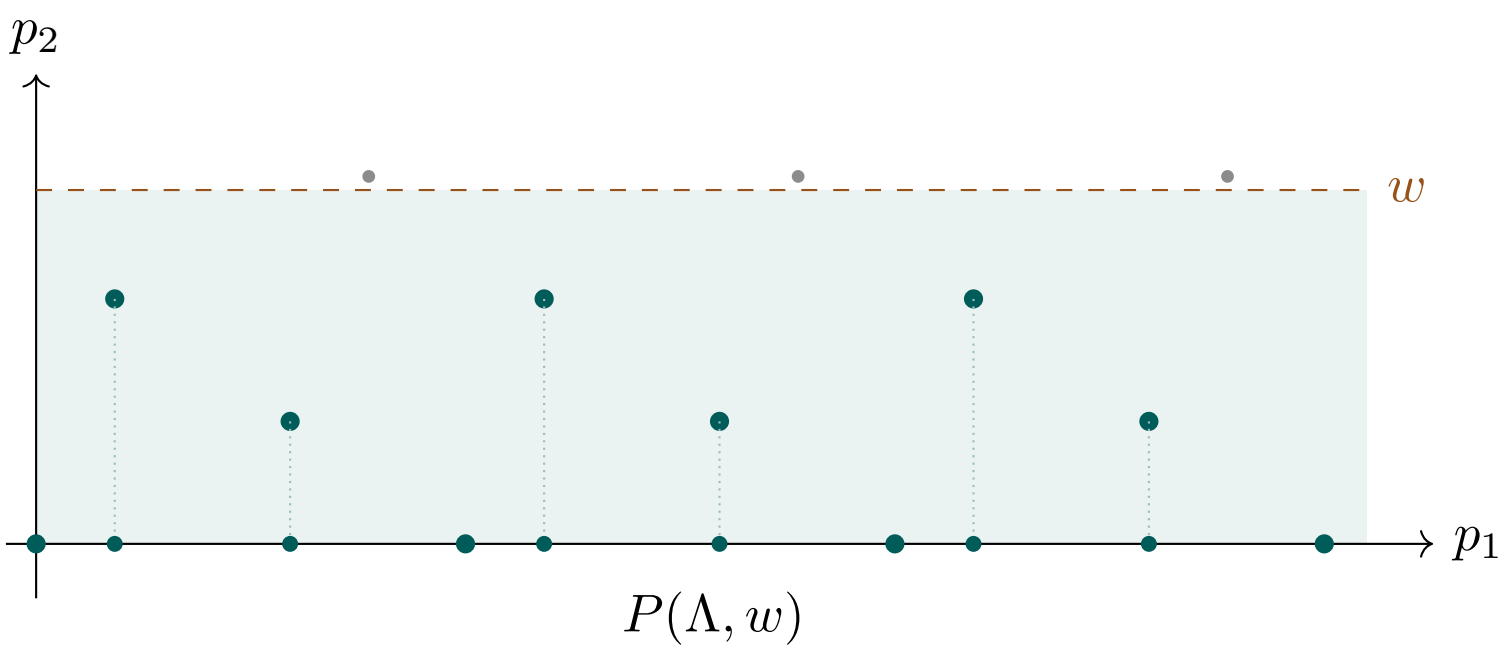
\begin{figure}[htbp]
\centering
\begin{tikzpicture}[x=.8cm,y=.9cm]
\fill[conn!8] (0,0) rectangle (11,2.6);
\draw[->] (-.25,0)--(11.55,0) node[right]{$p_1$};
\draw[->] (0,-.4)--(0,3.45) node[above]{$p_2$};
\draw[dashed,bound] (0,2.6)--(11,2.6) node[right]{$w$};
\begin{scope}
\clip (-.2,-.3) rectangle (11,3.15);
\foreach \i in {-3,...,10}{\foreach \j in {-3,...,8}{\pgfmathsetmacro{\xx}{1.45*\i+.65*\j}\pgfmathsetmacro{\yy}{-.9*\i+1.8*\j}\pgfmathtruncatemacro{\inside}{(\yy>=-.0001)&&(\yy<2.5999)}\ifnum\inside=1\relax
\fill[conn] (\xx,\yy) circle (1.8pt);\draw[conn!40,densely dotted] (\xx,\yy)--(\xx,0);\fill[conn] (\xx,0) circle (1.5pt);\else\fill[black!45] (\xx,\yy) circle (1.2pt);\fi}}
\end{scope}
\node[below] at (5.6,-.2) {$P(\Lambda,w)$};
\end{tikzpicture}
\caption{Projection of the lattice points in a strip. The projected points are ordered by their physical coordinates (schematic).}
\label{fig:strip}
\end{figure}

\subsection{The AF subgroupoid}
We recall the orbit-breaking construction (see \cite{Put89}). If $S$ is a minimal homeomorphism of a Cantor set $\mathcal{X}$ and $C$ is a nonempty clopen subset, then the first-return times to $C$ are finite and locally constant. Deleting the arrows of $S$ whose ranges belong to $C$ gives a compact open elementary subgroupoid of $\mathcal{X}\rtimes_S\mathbb Z$: its orbits consist of the finite intervals between consecutive returns to $C$. As the clopen sets $C_i$ decrease, these subgroupoids increase. If their intersection is a finite set $K$ meeting each orbit of $S$ at most once, then the union is generated by the arrows $(S,\zeta)$ with $S(\zeta)\notin K$. It is AF, and it is minimal since each positive or negative half-orbit is dense. For $m\ne0$, precisely $|m||K|$ arrows in the bisection of $S^m$ are missing from the union. Indeed, a path of $|m|$ consecutive steps crosses at most one cut, and each cut determines $|m|$ distinct initial points. Thus, the complement is discrete.

Taking $C_i$ to be a clopen neighborhood basis of the distinguished unit $\xi$ proves the first assertion of Proposition~\ref{main:reconstruction}. For reduced lattices, the AF subgroupoid defined by the clopen orbit transversals arising from modular returns need not coincide with this subgroupoid.

\subsection{Recovery of the lattice and window}
We give a direct calculation in the specified physical and internal coordinates (see \cite[Sections~4--5]{BLM07}). The cocycle at $\xi$ recovers $P=\{c(g):\src(g)=\xi\}$. Put $D=\langle P\rangle_{\mathbb Z}$, write $D_0=p_1(\Lambda)$, and define $\iota:D_0\to\mathbb R$ by $(d,\iota(d))\in\Lambda$.

\begin{lemma}\label{lem:cone}
The additive monoid $D^+$ generated by $P$ satisfies $D^+=\{d\in D_0:\iota(d)\geq0\}$. In particular, $D=D_0$.
\end{lemma}
\begin{proof}
The inclusion $D^+\subseteq\{d\in D_0:\iota(d)\geq0\}$ follows from the definition of $P$. Choose $z\in D_0$ with $0<\iota(z)<w$. For $d\in D_0$ with $\iota(d)\geq0$, take $m=\lfloor\iota(d)/\iota(z)\rfloor$. Then $z$ and $d-mz$ belong to $P$, and $d=mz+(d-mz)$. Applying the same argument to $-d$ for elements with negative internal coordinate gives $D=D_0$.
\end{proof}

The order defined by $D^+$ determines $\iota$ up to a positive scalar. Indeed, for a positive element $e$, comparisons of $nd$ with $me$ for integers $m$ and positive integers $n$ recover $\iota(d)/\iota(e)$. If $(d_1,d_2)$ is a basis of $D$, then the equation
\[
 |d_1\iota(d_2)-d_2\iota(d_1)|=1
\]
fixes the scalar. We recover $\Lambda=\{(d,\iota(d)):d\in D\}$ and $w=\sup_{p\in P}\iota(p)$, the latter by density of $p_2(\Lambda)$. This proves Proposition~\ref{main:reconstruction}. The choice of the half-open interval fixes which internal order is recovered; for an internal reflection one uses the reflected window and the corresponding side of its endpoints.

For the diagonal deformation $A_t$ described in \subsecref{subsec:prelim-modular}, we have
\begin{equation}\label{eq:flowstrip}
 P(A_t\Lambda,e^{-t/2}w)=e^{t/2}P(\Lambda,w).
\end{equation}
The induced conjugacy of the first return maps multiplies $c$ by $e^{t/2}$. Thus, the groupoid without the cocycle does not distinguish the diagonal deformation. See also \cite{CS03,CS06,KS17} for deformations of tilings and reprojections of model sets.

\section{Modular returns and restrictions of groupoids}\label{sec:returns}
\subsection{Reduced lattices and the translation flow}\label{subsec:reduced}
With the conventions of \subsecref{subsec:prelim-modular}, we use the coordinates $(x,y)$ of the natural extension of the regular Gauss map, with $0<x,y<1$ irrational. The two orientation sheets are indexed by $\varepsilon\in\{1,-1\}$. On each sheet, we orient the internal coordinate according to the convention for the strip. The basis and its determinant are
\begin{equation}\label{eq:basis}
 u=(h,-yk),\qquad v=(xh,k),\qquad hk(1+xy)=1,\qquad h,k>0.
\end{equation}
Reflection of the internal coordinate exchanges the sheets. We retain the sheet and the one-sided convention for the window of width $k$. On the reduced section, $h=k=(1+xy)^{-1/2}$; the formulas also apply to diagonal deformations. For the natural extension and its realization by cross-sections for the geodesic flow, see \cite{Ser85,Arn94,AS13,AS26}. The subdivision of the strip under the return map is described by the scenery flow for one-dimensional tilings \cite{AF01,AF05}. The basis \eqref{eq:basis} specifies the displacement cocycles, the distinguished units and the normalization of return times.

Let $t\in[0,k)$ be the internal coordinate of a point of the strip. If $t\geq yk$, then the next point of the strip to the right is obtained by adding $u$; otherwise, it is obtained by adding $u+v$. Indeed, any lattice vector to a closer point on the right has normalized coordinates $(m+nx,n-my)$ with $m,n\in\mathbb Z$, and satisfies
\[
 0<m+nx\leq1+x,\qquad |n-my|<1.
\]
If $n<0$, then $m>0$ and $n-my<-1$; if $m<0$, then $n>0$ and $n-my>1$. Thus, $m,n\geq0$. The first inequality rules out $m\geq2$. If $n\geq2$, then the second requires $m\geq2$. Hence, $m,n\in\{0,1\}$. The vector $v$ has internal displacement $k$ and is excluded. Exactly one of $u$ and $u+v$ stays in the window, with corresponding one-sided limits at the endpoints.

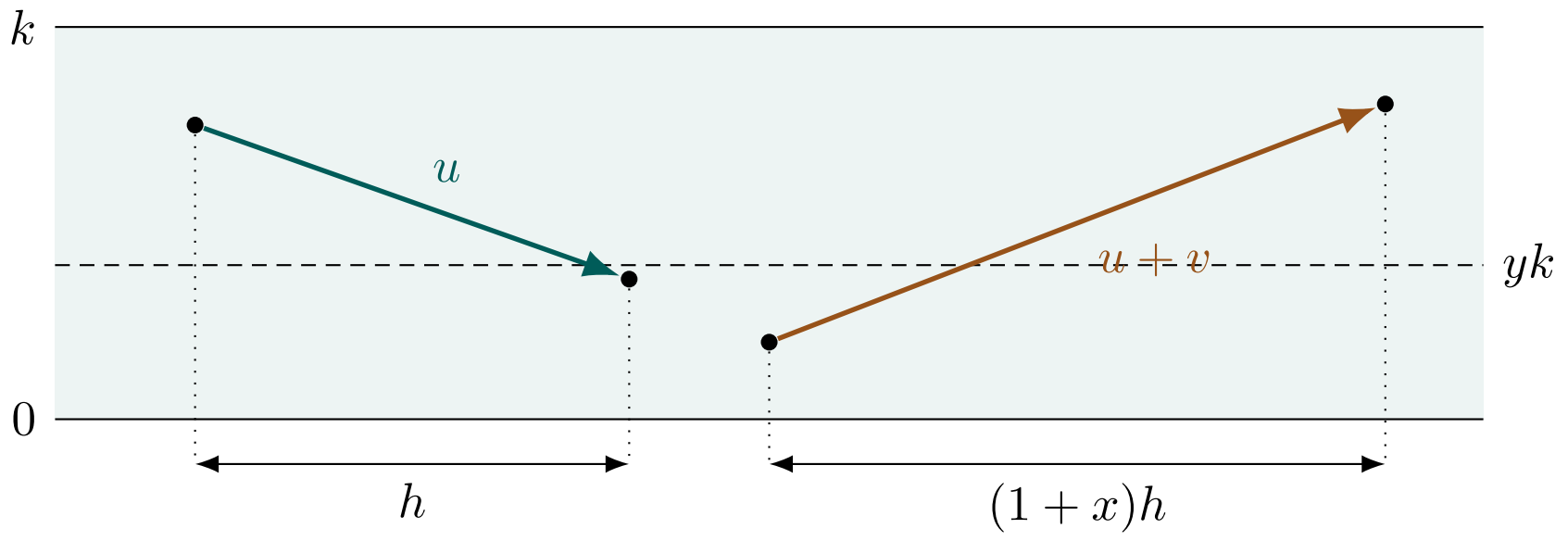
\begin{figure}[htbp]
\centering
\begin{tikzpicture}[x=1cm,y=1cm,>=Latex]
\fill[conn!7] (0,0) rectangle (10.2,2.8);
\draw (0,0)--(10.2,0); \draw (0,2.8)--(10.2,2.8);
\draw[densely dashed] (0,1.1)--(10.2,1.1) node[right]{$yk$};
\node[left] at (0,0){$0$};\node[left] at (0,2.8){$k$};
\node[point] (u1) at (1,2.1) {}; \node[point] (u2) at (4.1,1) {};
\draw[->,conn,line width=1pt] (u1)--(u2) node[midway,above right]{$u$};
\node[point] (v1) at (5.1,.55) {}; \node[point] (v2) at (9.5,2.25) {};
\draw[->,bound,line width=1pt] (v1)--(v2) node[midway,below right]{$u+v$};
\draw[dotted] (u1)--(1,-.32) (u2)--(4.1,-.32) (v1)--(5.1,-.32) (v2)--(9.5,-.32);
\draw[<->] (1,-.32)--(4.1,-.32) node[midway,below]{$h$};
\draw[<->] (5.1,-.32)--(9.5,-.32) node[midway,below]{$(1+x)h$};
\end{tikzpicture}
\caption{Vectors between consecutive points of a reduced strip. The displacements in the internal direction are $-yk$ and $(1-y)k$.}
\label{fig:successors}
\end{figure}

The action on the internal coordinate after normalization is given by $S_y(t)=t-y\pmod 1$. The intervals between consecutive projected points have two types, $\mathrm S$ and $\mathrm L$, with lengths $\ell_S=h$ and $\ell_L=(1+x)h$, respectively. The corresponding intervals of internal coordinates have lengths $1-y$ and $y$. Splitting the orbit $\mathbb Zy\pmod 1$ of the endpoints gives $\mathcal{X}_y$. The two-sided coding by $[0,y)$ and $[y,1)$ extends to $\mathcal{X}_y$ by one-sided limits. It is injective because the translates of the one-sided intervals separate points. By density of the orbit and the occurrence of both limits, its image is the closure of the translates of the gap sequence of the strip. Thus, we identify $\mathfrak G_{\Lambda,k}$ with $\mathfrak N_y=\mathcal{X}_y\rtimes_{S_y}\mathbb Z$. This is the Sturmian coding of the rotation (see \cite{MH40,GMP03}).

At each unit, the arrow corresponding to $S_y$ has the least positive displacement. The unique invariant probability measure $\mu_y$ of $\mathfrak N_y$ is described in \subsecref{subsec:prelim-frag}. If $U_L$ is the clopen set at which the first return time of the translation flow is $\ell_L$, then $\mu_y(U_L)=y$. Thus, the parameters in \eqref{eq:basis} are recovered without a distinguished unit by
\begin{equation}\label{eq:reducedrecovery}
 h=\ell_S,\qquad x=\frac{\ell_L}{\ell_S}-1,\qquad y=\mu_y(U_L),\qquad k=\frac{1}{h(1+xy)}.
\end{equation}
The distinguished unit in Proposition~\ref{main:reconstruction} specifies the pointed configuration; for modular histories, we also retain the sheet.

For example, if $|b|$ is sufficiently small, then the shear $(s,t)\mapsto(s+bt,t)$ preserves the coding intervals and changes the lengths to $h-byk$ and $(1+x)h+b(1-y)k$. Their weighted mean remains $h(1+xy)=1/k$. Thus, even the groupoid together with the mean return time of the translation flow need not determine the original lattice. The exact cocycle records both lengths.

\subsection{Induced transformations and the displacement cocycle}\label{subsec:induced}
Put $a_0=\lfloor1/x\rfloor$, $x'=T(x)=1/x-a_0$, and $y'=1/(a_0+y)$, where $T$ is the Gauss map of \subsecref{subsec:prelim-modular}. We use the regular case of the cross-section of Arnoux \cite{Arn94}, in the coordinates of \cite[Theorem~4.1 and Lemma~4.2]{AS13}. With the sheet retained, the return map $\widehat T$ projects to $T$ in the first coordinate:
\begin{equation}\label{eq:gauss}
 \widehat T(x,y,\varepsilon)=(x',y',-\varepsilon).
\end{equation}
The change of basis from $(u,v)$ to $(v,u-a_0v)$ reverses orientation. Reversing the internal coordinate in the new sheet gives the basis
\[
 \widetilde u=(xh,-k),\qquad \widetilde v=(xh x',(a_0+y)k).
\]
Before applying the diagonal deformation, we have $\widetilde h=xh$ and $\widetilde k=(a_0+y)k$. In the section $h=k$, the flow $A_r$ makes these quantities equal at
\begin{equation}\label{eq:roof}
 r_{\rm lat}(x,y)=-\log(xy').
\end{equation}
The return time in the normalization of \cite[Lemma~5.2]{AS13} is $\tau(x)=-2\log x$. With $L(x,y)=\log(1+xy)$, the identities $1+x'y'=(1+xy)y'/x$ and $-\log y'=-\log x+L(x,y)-L(x',y')$ give
\begin{equation}\label{eq:coboundary}
 r_{\rm lat}(x,y)=\tau(x)+L(x,y)-L(x',y').
\end{equation}
Since $0<L<\log2$, the accumulated return times in the two normalizations differ by a bounded amount. With the convention $(z,t+r_{\rm lat}(z))\sim(\widehat Tz,t)$, the map $(z,t)\mapsto(z,t-L(z))$ conjugates the suspension with roof $r_{\rm lat}$ to that with roof $\tau$. We exclude rational boundary configurations.

\begin{figure}[htbp]
\centering
\begin{tikzpicture}[x=2.12cm,y=2.12cm,>=Latex]
\draw[thin] (-3.8,0)--(1.58,0);
\draw[black!45,thin] (-3.5,0.8660254) arc[start angle=120,end angle=60,radius=1];
\draw[black!45,thin] (-3.5,0.8660254)--(-3.5,2.24);
\draw[black!45,thin] (-2.5,0.8660254) arc[start angle=120,end angle=60,radius=1];
\draw[black!45,thin] (-2.5,0.8660254)--(-2.5,2.24);
\draw[black!45,thin] (-1.5,0.8660254) arc[start angle=120,end angle=60,radius=1];
\draw[black!45,thin] (-1.5,0.8660254)--(-1.5,2.24);
\draw[black!45,thin] (-0.5,0.8660254) arc[start angle=120,end angle=60,radius=1];
\draw[black!45,thin] (-0.5,0.8660254)--(-0.5,2.24);
\draw[black!45,thin] (0.5,0.8660254) arc[start angle=120,end angle=60,radius=1];
\draw[black!45,thin] (0.5,0.8660254)--(0.5,2.24);
\draw[black!45,thin] (1.5,0.8660254)--(1.5,2.24);
\draw[black!50] (-3,0)--(-3,-.035);
\node[below,font=\scriptsize] at (-3,-.02) {$-3$};
\draw[black!50] (-2,0)--(-2,-.035);
\node[below,font=\scriptsize] at (-2,-.02) {$-2$};
\draw[black!50] (-1,0)--(-1,-.035);
\node[below,font=\scriptsize] at (-1,-.02) {$-1$};
\draw[black!50] (0,0)--(0,-.035);
\node[below,font=\scriptsize] at (0,-.02) {$0$};
\draw[black!50] (1,0)--(1,-.035);
\node[below,font=\scriptsize] at (1,-.02) {$1$};
\draw[black!70,line width=.55pt] (0.6180340,0.0000000) arc[start angle=0.0000000,end angle=180.0000000,radius=2.0161238];
\draw[line width=1.05pt,-{Latex[length=2mm]}] (0.1871161,1.2457437) arc[start angle=38.1622497,end angle=95.1527480,radius=2.0161238];
\draw[line width=1.05pt] (-1.5791604,2.0079762) arc[start angle=95.1527480,end angle=152.1432463,radius=2.0161238];
\node[point] at (-3.4142136,0) {};
\node[below,font=\small] at (-3.4142136,-.02) {$-1/x$};
\node[point] at (0.6180340,0) {};
\node[below,font=\small] at (0.6180340,-.02) {$y$};
\node[point] at (0.1871161,1.2457437) {};
\draw[thin] (0.1871161,1.2457437)--(0.46,1.44) node[right,font=\scriptsize]{$M_-i$};
\node[point] at (-3.1805823,0.9420593) {};
\draw[thin] (-3.1805823,0.9420593)--(-3.38,1.36) node[left,font=\scriptsize]{$M_-e^{\tau(x)}i$};
\node[above,font=\small] at (-1.1099999999999999,2.2800000000000002) {$\varepsilon=-1$};
\end{tikzpicture}
\caption{A lift of the geodesic from $y$ to $-1/x$ on the first sheet. The thick arc is traversed during the return time $\tau(x)$. Here $x=1/(2+\sqrt2)$, $y=(\sqrt5-1)/2$, $a_0=3$, and $M_-=M_-(x,y)$. The thin gray curves are boundaries of translates of the usual modular fundamental domain.}
\label{fig:gaussgeodesic}
\end{figure}
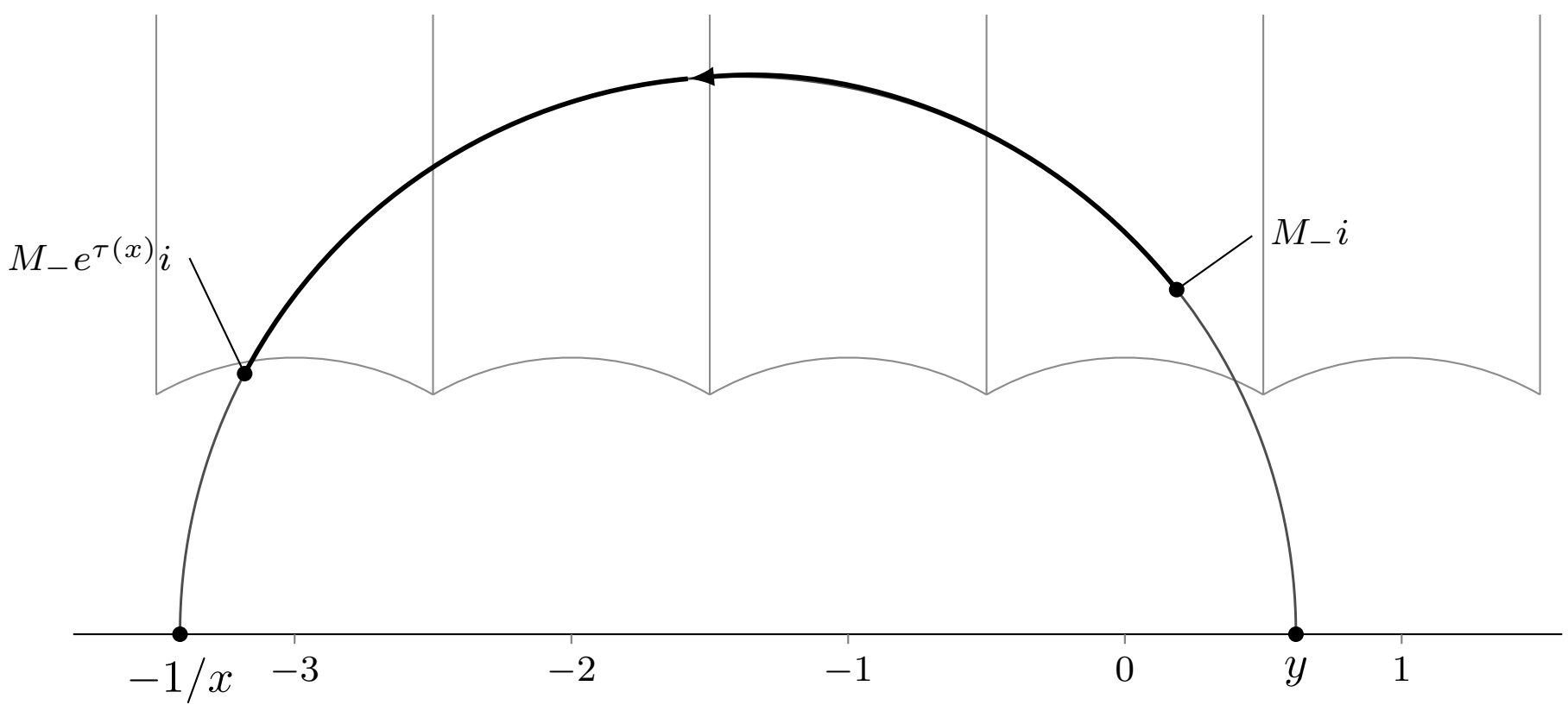

To illustrate the return map, we use the upper half-plane representatives of \cite[Subsection~4.2 and Lemma~5.2]{AS13}. With $\beta=y/(1+xy)$, these representatives are the matrices
\[
 M_-(x,y)=\begin{pmatrix}1&\beta\\-x&1-x\beta\end{pmatrix},\qquad
 M_+(x,y)=\begin{pmatrix}x&1-x\beta\\-1&\beta\end{pmatrix}.
\]
On the sheet $\varepsilon=-1$, the geodesic has initial and terminal endpoints $y$ and $-1/x$, respectively. During one return, the geodesic goes from $M_-(x,y)i$ to $M_-(x,y)(e^{\tau(x)}i)$ (Figure~\ref{fig:gaussgeodesic}). Translation by $a_0$ sends the endpoints to $a_0+y=1/y'$ and $-x'$, and the terminal point of the segment to $M_+(x',y')i$.

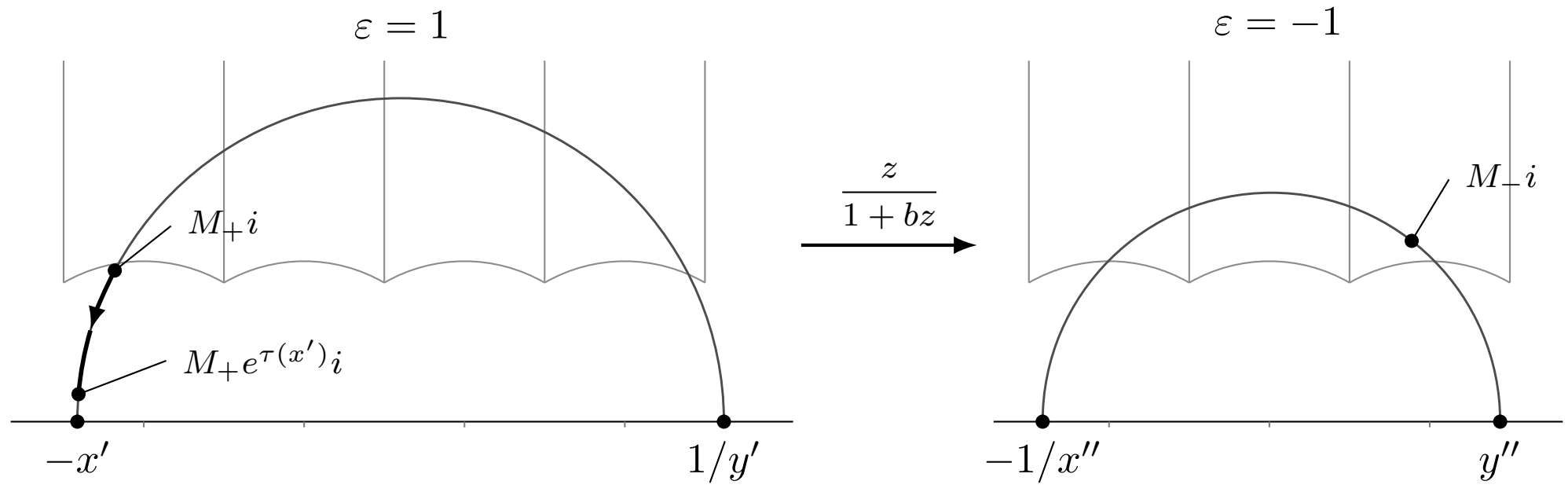
\begin{figure}[htbp]
\centering
\begin{tikzpicture}[x=1.35cm,y=1.35cm,>=Latex]
\draw[thin] (-0.83,0)--(4.05,0);
\draw[black!45,thin] (-0.5,0.8660254) arc[start angle=120,end angle=60,radius=1];
\draw[black!45,thin] (-0.5,0.8660254)--(-0.5,2.25);
\draw[black!45,thin] (0.5,0.8660254) arc[start angle=120,end angle=60,radius=1];
\draw[black!45,thin] (0.5,0.8660254)--(0.5,2.25);
\draw[black!45,thin] (1.5,0.8660254) arc[start angle=120,end angle=60,radius=1];
\draw[black!45,thin] (1.5,0.8660254)--(1.5,2.25);
\draw[black!45,thin] (2.5,0.8660254) arc[start angle=120,end angle=60,radius=1];
\draw[black!45,thin] (2.5,0.8660254)--(2.5,2.25);
\draw[black!45,thin] (3.5,0.8660254)--(3.5,2.25);
\draw[black!50] (0,0)--(0,-.035);
\draw[black!50] (1,0)--(1,-.035);
\draw[black!50] (2,0)--(2,-.035);
\draw[black!50] (3,0)--(3,-.035);
\draw[black!70,line width=.55pt] (3.6180340,0.0000000) arc[start angle=0.0000000,end angle=180.0000000,radius=2.0161238];
\draw[line width=1.05pt,-{Latex[length=2mm]}] (-0.1805823,0.9420593) arc[start angle=152.1432463,end angle=163.6351469,radius=2.0161238];
\draw[line width=1.05pt] (-0.3325343,0.5680488) arc[start angle=163.6351469,end angle=175.1270475,radius=2.0161238];
\node[point] at (-0.4142136,0) {};
\node[below,font=\small] at (-0.4142136,-.02) {$-x'$};
\node[point] at (3.6180340,0) {};
\node[below,font=\small] at (3.6180340,-.02) {$1/y'$};
\node[point] at (-0.1805823,0.9420593) {};
\draw[thin] (-0.1805823,0.9420593)--(0.17,1.22) node[right,font=\scriptsize]{$M_+i$};
\node[point] at (-0.4069263,0.1712628) {};
\draw[thin] (-0.4069263,0.1712628)--(0.14,0.38) node[right,font=\scriptsize]{$M_+e^{\tau(x')}i$};
\node[above,font=\small] at (1.6099999999999999,2.29) {$\varepsilon=1$};
\draw[-{Latex[length=2mm]},thick] (4.10,1.10)--(5.20,1.10);
\node[above,font=\scriptsize] at (4.65,1.10){$\displaystyle\frac{z}{1+bz}$};
\begin{scope}[shift={(8.02,0)}]
\draw[thin] (-2.72,0)--(0.83,0);
\draw[black!45,thin] (-2.5,0.8660254) arc[start angle=120,end angle=60,radius=1];
\draw[black!45,thin] (-2.5,0.8660254)--(-2.5,2.25);
\draw[black!45,thin] (-1.5,0.8660254) arc[start angle=120,end angle=60,radius=1];
\draw[black!45,thin] (-1.5,0.8660254)--(-1.5,2.25);
\draw[black!45,thin] (-0.5,0.8660254) arc[start angle=120,end angle=60,radius=1];
\draw[black!45,thin] (-0.5,0.8660254)--(-0.5,2.25);
\draw[black!45,thin] (0.5,0.8660254)--(0.5,2.25);
\draw[black!50] (-2,0)--(-2,-.035);
\draw[black!50] (-1,0)--(-1,-.035);
\draw[black!50] (0,0)--(0,-.035);
\draw[black!70,line width=.55pt] (0.4392914,0.0000000) arc[start angle=0.0000000,end angle=180.0000000,radius=1.4267525];
\node[point] at (-2.4142136,0) {};
\node[below,font=\small] at (-2.4142136,-.02) {$-1/x''$};
\node[point] at (0.4392914,0) {};
\node[below,font=\small] at (0.4392914,-.02) {$y''$};
\node[point] at (-0.1124291,1.1269169) {};
\draw[thin] (-0.1124291,1.1269169)--(0.12,1.51) node[right,font=\scriptsize]{$M_-i$};
\node[above,font=\small] at (-0.9450000000000001,2.29) {$\varepsilon=-1$};
\end{scope}
\end{tikzpicture}
\caption{After translation by $a_0$, the geodesic in Figure~\ref{fig:gaussgeodesic} is represented on the other sheet (left). Its next return, of duration $\tau(x')$, is brought to the first sheet by $z\mapsto z/(1+bz)$ (right), where $b=\lfloor1/x'\rfloor=2$. In the two panels $M_+=M_+(x',y')$ and $M_-=M_-(x'',y'')$, respectively.}
\label{fig:gaussnext}
\end{figure}

For the next return, put $b=\lfloor1/x'\rfloor$, $x''=1/x'-b$, and $y''=1/(b+y')$. The modular transformation $z\mapsto z/(1+bz)$ sends the endpoints $1/y'$ and $-x'$ to $y''$ and $-1/x''$, respectively, and sends $M_+(x',y')(e^{\tau(x')}i)$ to $M_-(x'',y'')i$ (Figure~\ref{fig:gaussnext}). The points marked in both figures correspond to the return times $\tau$ rather than $r_{\rm lat}$ (see \eqref{eq:coboundary}).

In $\mathcal{X}_{y'}$, take the clopen interval $C=[0^+,(y')^-]$. Let $F$ be the induced transformation of $S_{y'}$ on $C$. For $0<t<y'$, if $t<(1-y)y'$, then the first-return path has $a_0$ steps; if $t>(1-y)y'$, then it has $a_0+1$ steps. Indeed, the first step adds $1-y'$, and the subsequent steps subtract $y'$ until the path returns. The two formulas for $F$ are $t\mapsto t+yy'$ and $t\mapsto t-(1-y)y'$, respectively. They extend to the corresponding one-sided lifts. The affine map
\begin{equation}\label{eq:phi}
 \phi:\mathcal{X}_y\longrightarrow C,\qquad \phi(t)=y'(1-t),
\end{equation}
where the two lifts of $0$ are represented by the endpoints of $[0,1]$, is a conjugacy from $S_y$ to $F$. The equality $1/y'=a_0+y$ shows that $t\in\mathbb Z+y\mathbb Z$ if and only if $y'(1-t)\in\mathbb Z+y'\mathbb Z$. The map interchanges the two sides at split points. The distinguished point $0^+$ is sent to $(y')^-$, while $\phi(0^-)=0^+$.

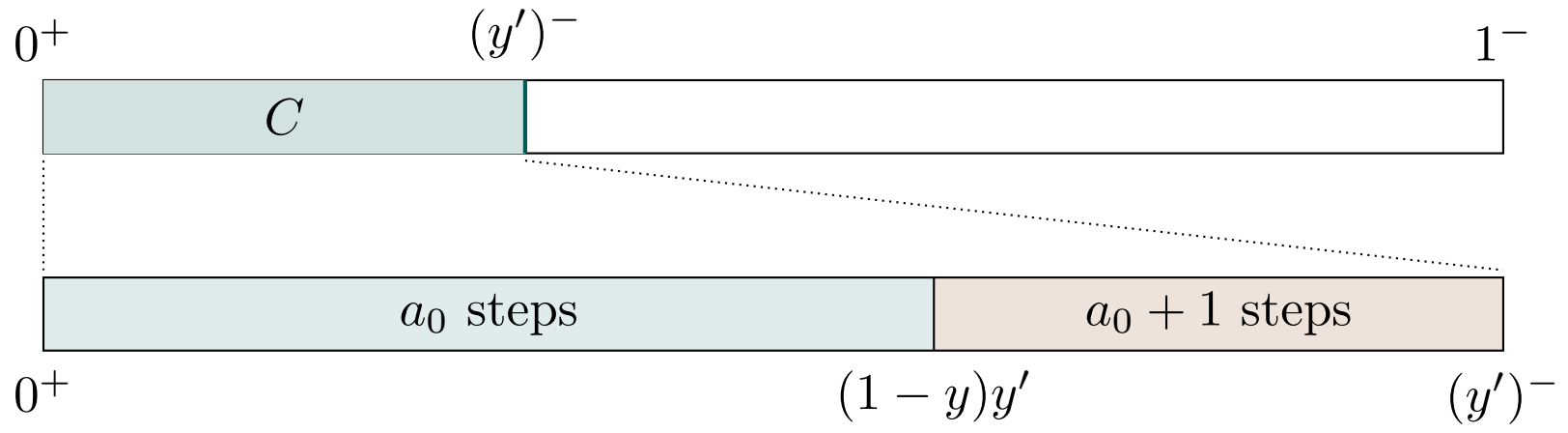
\begin{figure}[htbp]
\centering
\begin{tikzpicture}[x=1cm,y=1cm,>=Latex]
\draw (0,1.7) rectangle (10,2.2);\fill[conn!18] (0,1.7) rectangle (3.3,2.2);\draw[conn,thick] (3.3,1.7)--(3.3,2.2);
\node[above] at (0,2.2){$0^+$};\node[above] at (3.3,2.2){$(y')^-$};\node[above] at (10,2.2){$1^-$};
\node at (1.65,1.95){$C$};
\draw[densely dotted] (0,1.65)--(0,.9) (3.3,1.65)--(10,.9);
\fill[conn!12] (0,.35) rectangle (6.1,.85);\fill[bound!16] (6.1,.35) rectangle (10,.85);\draw (0,.35) rectangle (10,.85);\draw (6.1,.35)--(6.1,.85);
\node at (3.05,.6){$a_0$ steps};\node at (8.05,.6){$a_0+1$ steps};
\node[below] at (0,.35){$0^+$};\node[below] at (6.1,.35){$(1-y)y'$};\node[below] at (10,.35){$(y')^-$};
\end{tikzpicture}
\caption{The inducing interval of $S_{y'}$ and its two first-return times. The drawing of the interval is schematic.}
\label{fig:induction}
\end{figure}

As in the tiling description of the scenery flow \cite{AF01,AF05}, before applying the diagonal deformation, the two new physical lengths are $xh$ and $(1+x')xh$. Each first-return path starts with a long step and then takes only short steps. With the new types denoted by $\mathrm S'$ and $\mathrm L'$, under \eqref{eq:phi} the old tiles correspond to the words
\begin{equation}\label{eq:inflation}
 \mathrm S\longmapsto\mathrm L'(\mathrm S')^{a_0-1},\qquad \mathrm L\longmapsto\mathrm L'(\mathrm S')^{a_0}.
\end{equation}
Their physical displacements are $(a_0+x')xh=h$ and $(a_0+1+x')xh=(1+x)h$, respectively. This also holds for $a_0=1$. The induced cocycle agrees with the original displacement cocycle before rescaling and is multiplied by $e^{r_{\rm lat}/2}$ after rescaling.

\begin{figure}[htbp]
\centering
\begin{tikzpicture}[x=.9cm,y=.8cm]
\node[left] at (0,1.5){$S$};\node[left] at (0,0){$L$};
\foreach \y/\n in {1.5/2,0/3}{
\fill[bound!18] (0,\y-.36) rectangle (2.4,\y+.36);\draw (0,\y-.36) rectangle (2.4,\y+.36);\node at (1.2,\y){$L'$};
\foreach \i in {1,...,\n}{\pgfmathsetmacro{\a}{2.4+1.6*(\i-1)}\pgfmathsetmacro{\b}{\a+1.6}\fill[conn!10] (\a,\y-.36) rectangle (\b,\y+.36);\draw (\a,\y-.36) rectangle (\b,\y+.36);\node at ({(\a+\b)/2},\y){$S'$};}}
\end{tikzpicture}
\caption{The subdivisions of the two types of intervals between consecutive projected points for $a_0=3$.}
\label{fig:ordered}
\end{figure}
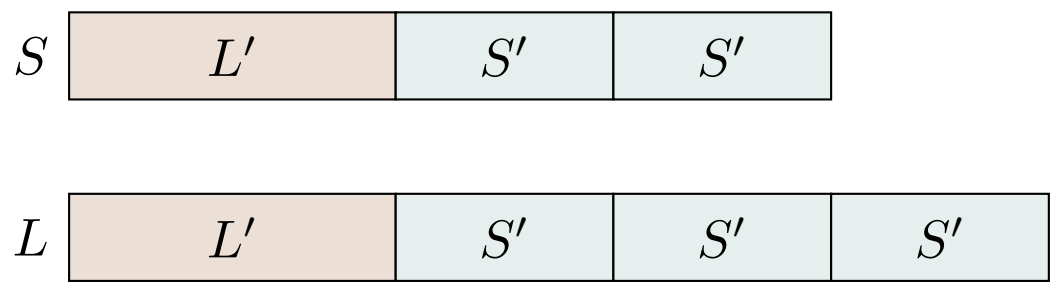

The conjugacy \eqref{eq:phi} induces an isomorphism $\Phi$ onto the restriction of the new groupoid to $C$, since every arrow with source and range in $C$ is given by an iterate of the first-return transformation. Since $S_{y'}$ is minimal, $C$ is an orbit transversal, and hence the groupoids are Morita equivalent. For the cocycles and the distinguished units, we have
\begin{equation}\label{eq:markedind}
 \Phi:\mathfrak N_y\xrightarrow{\ \cong\ }\mathfrak N_{y'}|_C,\qquad
 c'(\Phi(g))=e^{r_{\rm lat}/2}c(g),\qquad \Phi(\xi)=\phi(\xi).
\end{equation}
Here, $c'$ is the cocycle after rescaling the lattice. The factor depends on the normalization. Restriction to $C$ changes the ambient groupoid, whereas telescoping changes only a presentation of its tail groupoid.

The successive maximal blocks of equal symbols in a regular Farey cutting sequence have lengths given by the partial quotients (see \cite{Ser85,AS26}). For a block of length $a$, the matrix of the change of basis and that of the ordered decomposition \eqref{eq:inflation} are related by
\[
 B_a=\begin{pmatrix}0&1\\1&-a\end{pmatrix},\qquad
 E_a=\begin{pmatrix}a-1&a\\1&1\end{pmatrix}=H^{-1}B_a^{-1}H,
 \qquad H=\begin{pmatrix}1&1\\0&1\end{pmatrix}.
\]
The matrix $H$ changes the basis $(u,v)$ to $(u,u+v)$. The rows and columns of $E_a$ are indexed by $(\mathrm S',\mathrm L')$ and $(\mathrm S,\mathrm L)$, respectively. Thus, such a block determines the isomorphism \eqref{eq:markedind}, not an arrow of $\mathfrak N_y$ corresponding to one application of $S_y$.

\subsection{Successive modular returns}\label{subsec:compression}
Let $z_n=(x_n,y_n,\varepsilon_n)$, $n\in\mathbb Z$, be a two-sided orbit of \eqref{eq:gauss} with both continued-fraction coordinates irrational. Write $\phi_n$ for the map \eqref{eq:phi} between successive unit spaces. The composition is
\[
 \phi_{N-1}\circ\cdots\circ\phi_0(t)=d_N+(-1)^N\omega_Nt,
 \qquad \omega_N=\prod_{i=1}^N y_i,
\]
where $d_0=0$ and $d_{N+1}=y_{N+1}(1-d_N)$. Let $C_N$ be the image of $\mathcal{X}_{y_0}$ in $\mathcal{X}_{y_N}$. Composition gives $\Phi_N:\mathfrak N_{y_0}\to\mathfrak N_{y_N}|_{C_N}$, with $\mu_{y_N}(C_N)=\omega_N$. If $N$ is even, then the image of $0^+$ is the lower one-sided endpoint; otherwise, it is the upper one. The two lifts of $0$ lie in distinct $S_y$-orbits, so replacing the image by $0^+$ need not be realized by an arrow of the groupoid.

We express the classical continuant identities (see \cite{Khi64}) in terms of the restrictions.
\begin{proposition}\label{prop:compression}
If $p_N/q_N$ is the $N$-th convergent of $x_0=[0;a_1,a_2,\ldots]$, then
\begin{equation}\label{eq:exactcompression}
 \omega_N^{-1}=q_N+p_Ny_0.
\end{equation}
Moreover, if $s_N=\exp(\frac12\sum_{n=0}^{N-1}r_{{\rm lat},n})$, then
\begin{equation}\label{eq:scaleidentity}
 c_N(\Phi_N(g))=s_Nc_0(g),\qquad
 s_N\omega_N=\sqrt{\frac{1+x_Ny_N}{1+x_0y_0}}.
\end{equation}
In particular,
\begin{equation}\label{eq:compressroof}
 \left|-\log\mu_{y_N}(C_N)-\frac12\sum_{n=0}^{N-1}\tau_n\right|\leq\log 2.
\end{equation}
\end{proposition}
\begin{proof}
The numbers $K_N=\omega_N^{-1}$ satisfy $K_0=1$, $K_1=a_1+y_0$, and $K_{N+1}=a_{N+1}K_N+K_{N-1}$ by the recursion for $y_n$. This is the recurrence for $q_N+p_Ny_0$, proving \eqref{eq:exactcompression}. The first identity of \eqref{eq:scaleidentity} follows by composition of \eqref{eq:markedind}. For the second, sum the identity $-\log y_{n+1}=-\log x_n+L(x_n,y_n)-L(x_{n+1},y_{n+1})$ and use \eqref{eq:roof}. The bound in \eqref{eq:compressroof} follows from the same sum and $0<L<\log2$.
\end{proof}

Composing the embeddings from $-N$ to $0$ gives nested intervals $C_N^-\subseteq\mathcal{X}_{y_0}$ with
\begin{equation}\label{eq:backwardcompression}
 \mu_{y_0}(C_N^-)=\prod_{i=-N+1}^{0}y_i
 =\frac1{q_N(y_0)+q_{N-1}(y_0)T^N(y_0)},
\end{equation}
where the last equality is the continued-fraction product identity, since $y_{-i}=T^i(y_0)$. The projected intervals shrink to a point, so their intersection $K$ consists of one or two points of the split circle. In the latter case, these are the two lifts of a split point and lie in distinct $S_{y_0}$-orbits. The orbit-breaking construction in \secref{sec:strip} gives a minimal AF subgroupoid of $\mathfrak N_{y_0}$ with discrete complement by deleting the arrows of $S_{y_0}$ entering $K$. For instance, if $y=(\sqrt5-1)/2$, then the affine maps $t\mapsto y(1-t)$ have fixed point $1-y\in\mathbb Zy\pmod 1$ in the interior of every projected interval, so $K$ consists of both lifts of this point.

Proposition~\ref{main:reconstruction} and \eqref{eq:reducedrecovery} recover the reduced lattices. The isomorphisms \eqref{eq:markedind} record the ordered return paths, images of the distinguished units, and factors multiplying the cocycles. Together with the orientation sheet and suspension coordinate, these data recover the regular modular history by \eqref{eq:gauss} and its inverse. Conversely, the history determines the marked groupoids and these isomorphisms.

\subsection{Application to the growth of groups}\label{sec:modulargrowth}
Write $(x_n,y_n)_{n\in\mathbb Z}$ for a two-sided orbit of the natural extension with irrational coordinates and put $\tau_n=-2\log x_n$. Let $C_N^-$ be the clopen orbit transversals described in \subsecref{subsec:compression}. We apply Theorem~\ref{main:growth} to the fragmentation in \secref{sec:construction}.

\begin{corollary}\label{cor:modular}
If this two-sided orbit satisfies
\[
 \liminf_{N\to\infty}(2/3)^{N/64}\sum_{n=-N}^{-1}\tau_n=0,
\]
or, equivalently, $\liminf_{N\to\infty}(2/3)^{N/64}(-\log\mu_{y_0}(C_N^-))=0$, then $G_{y_0}$, obtained from the action groupoid of the first return map for the horizontal strip, has intermediate growth. If the averages of the backward return times are bounded, then the upper bound in Theorem~\ref{main:growth} holds with
\[
 \Lambda_{y_0}=\limsup_{N\to\infty}\frac{-\log\mu_{y_0}(C_N^-)}{N}.
\]
For the strip obtained by exchanging the coordinate directions with a sign change, the corresponding assertions hold for the forward return times and the group $G_{x_0}$, with $\Lambda_{x_0}=\limsup_{N\to\infty}(2N)^{-1}\sum_{n=0}^{N-1}\tau_n$.
\end{corollary}

By \subsecref{subsec:reduced}, the horizontal strip at time zero has rotation parameter $y_0$. The backward coordinates satisfy $y_{-n}=T^n(y_0)$. By \eqref{eq:backwardcompression},
\begin{equation}\label{eq:backlog}
 \log q_N(y_0)=-\log\mu_{y_0}(C_N^-)+O(1).
\end{equation}
Applying \eqref{eq:compressroof} to the returns from $-N$ to $0$ gives
\begin{equation}\label{eq:backroof}
 \log q_N(y_0)=\frac12\sum_{n=-N}^{-1}\tau_n+O(1).
\end{equation}
By \eqref{eq:coboundary}, both conditions are unchanged when we use the return times in the normalization $h=k$. Equations~\eqref{eq:backlog} and~\eqref{eq:backroof}, together with Theorem~\ref{main:growth}, prove the first assertion of Corollary~\ref{cor:modular}.

In the coordinates of \eqref{eq:basis}, let $J(s,t)=(t,-s)$. For $J\Lambda$, the basis $(Jv,-Ju)$ is
\[
 (k,-xh),\quad (yk,h),
\]
with the two parameters exchanged. Taking the width $h$, the strip in $J\Lambda$ has rotation parameter $x_0$. For the forward coordinates $x_n=T^n(x_0)$, the continued-fraction product identity gives
\begin{equation}\label{eq:forwardroof}
 \log q_N(x_0)=\frac12\sum_{n=0}^{N-1}\tau_n+O(1).
\end{equation}
This proves the second assertion of Corollary~\ref{cor:modular}. The conditions on forward and backward return times need not hold simultaneously.

\section{The fragmentation}\label{sec:construction}
We use the split circle $\mathcal{X}_\theta$ and the reflections $a,b$ of \subsecref{subsec:prelim-frag}, with the notation for convergents from \subsecref{subsec:prelim-modular}. Write $S_\theta$ for the rotation by $-\theta$ on $\mathcal{X}_\theta$. Its action groupoid $\mathfrak N_\theta=\mathcal{X}_\theta\rtimes_{S_\theta}\mathbb Z$, called the \emph{rotation groupoid}, is isomorphic to the action groupoid of the subgroup $\langle ba\rangle$. By \subsecref{subsec:reduced}, this is the action groupoid of the first return map for a reduced horizontal strip with $y=\theta$. The point $\eta=1/2$ is the unique fixed point of $a$ in $\mathcal{X}_\theta$. The two fixed points of $b$ are $\theta/2$ and $(1+\theta)/2$. All three are unsplit points, and their orbits under the dihedral group action are distinct. The dihedral group action is topologically free: every nonidentity rotation is fixed-point-free, and each reflection has at most two fixed points. Thus, its groupoid of germs coincides with its action groupoid. We denote this groupoid by $\mathfrak D_\theta$ and call it the \emph{dihedral groupoid}. The reflection induces an involutive automorphism $(r_\theta^m,\zeta)\mapsto(r_\theta^{-m},a(\zeta))$ of the rotation groupoid, where $r_\theta=ba$. The relation $ar_\theta a=r_\theta^{-1}$ gives
\[
 \mathfrak D_\theta\cong\mathfrak N_\theta\rtimes \mathbb Z/2\mathbb Z.
\]
In this presentation, the rotation groupoid is the kernel of the orientation cocycle. A reflection $x\mapsto p-x$ preserving $\mathbb Z\theta\pmod1$ has $p\in\mathbb Z\theta\pmod1$, and therefore gives the same dihedral groupoid. We use the reflection $a$ above to define the fragmentation.

\subsection{Choice of the denominators of convergents}\label{subsec:ranks}
For every $\theta$, let $n_1$ be the least index such that $q_{n_1}>200$, and set
\[
 n_j=n_1+32(j-1),\qquad Q_j=q_{n_j},\qquad D_j=q_{n_j+12}.
\]
The lower bounds on the denominators given by the Fibonacci sequence yield
\begin{equation}\label{eq:fib}
 D_j\geq233Q_j,\qquad Q_{j+1}\geq10946D_j.
\end{equation}
The bound on $Q_1$ places the intervals below in $(0,1)$, while the choices of $12$ and $32$ ensure that $200Q_j<D_j$ and $50/Q_{j+1}<1/(200D_j)$. These inequalities suffice for the separation estimates and the identification of the segments in \secref{sec:contraction}.
We use estimates for convergents (see \cite[Chapters~I--II]{Khi64}) and gaps between consecutive points (see \cite{Sla67}). By the next lemma, every block of $q_n$ consecutive iterates of any initial point meets every open interval of length greater than $3/q_n$. The same bound holds for clopen lifts, since they contain the preimages of the open intervals.

\begin{lemma}\label{lem:arithmetic}
For every denominator $q_n$ of a convergent, the gaps on the circle between the points of a block of $q_n$ consecutive iterates of the rotation have length at most $3/q_n$. If $q_n>1$, then the distance between any two distinct points of this block is greater than $1/(2q_n)$. Moreover,
\begin{align}
 \norm{k\theta}&>1/(2D_j)&& (0<|k|\leq200Q_j),\label{eq:sep1}\\
 \norm{\eta+k\theta}&>1/(4D_j)&& (|k|\leq100Q_j).\label{eq:sep2}
\end{align}
\end{lemma}
\begin{proof}
The points $k\theta$, $0\leq k<q_n$, are within $1/q_n$ of the equally spaced points $kp_n/q_n$. The lower bound on the distances follows from best approximation and $\norm{q_{n-1}\theta}>1/(q_n+q_{n-1})$. The same inequality for best approximation at the denominator $D_j$ gives \eqref{eq:sep1}, by \eqref{eq:fib}. For nonzero $k$ in \eqref{eq:sep2}, use $2\norm{\eta+k\theta}\geq\norm{2k\theta}$; the case $k=0$ is immediate.
\end{proof}

\subsection{Fragmentation of the reflection $a$}
For each $j$, choose four points $\eta+t_{j,i}$ of the orbit of $0$ such that
\[
 Q_jt_{j,i}\in [10+12i,14+12i],\qquad i=0,1,2,3.
\]
For each $i$, take the least $k\in\{0,\ldots,Q_j-1\}$ such that $k\theta\pmod1$ satisfies this condition for $\eta+t_{j,i}$. Such a $k$ exists by Lemma~\ref{lem:arithmetic}. The three intervening intervals have lengths at least $8/Q_j$. For $i=0,1,2$, let $P_{j,i}$ be the union of the clopen lift of the $i$-th interval (with the conventions of \subsecref{subsec:prelim-frag}) and its image under $a$. These sets are disjoint, also for distinct $j$ by \eqref{eq:fib}. We call $i$ the color of the $a$-edges with endpoints in $P_{j,i}$.

\begin{figure}[htbp]
\centering
\begin{tikzpicture}[x=.94cm,y=.78cm]
\foreach \row/\fac in {0/1,1/.46,2/.2}{
 \pgfmathsetmacro{\Y}{2.5-1.05*\row}
 \draw[black!45] (0,\Y)--(11,\Y);
 \node[left] at (0,\Y){$j\ifnum\row>0+\row\fi$};
 \foreach \i/\col in {0/conn,1/bound,2/blue!60!black}{
  \pgfmathsetmacro{\lo}{.9*\fac+1.35*\i*\fac}
  \pgfmathsetmacro{\hi}{\lo+1.35*\fac}
  \fill[\col!34] ({5.5+\lo},\Y-.2) rectangle ({5.5+\hi},\Y+.2);
  \fill[\col!34] ({5.5-\hi},\Y-.2) rectangle ({5.5-\lo},\Y+.2);
  \draw[\col] ({5.5+\lo},\Y-.2) rectangle ({5.5+\hi},\Y+.2);
  \draw[\col] ({5.5-\hi},\Y-.2) rectangle ({5.5-\lo},\Y+.2);
  \ifnum\row=0\relax\node at ({5.5+(\lo+\hi)/2},\Y){$\i$};\node at ({5.5-(\lo+\hi)/2},\Y){$\i$};\fi
 }
 \node[point] at (5.5,\Y) {};
}
\node[above] at (5.5,2.7){$\eta$};
\draw[densely dashed] (5.5,2.5)--(5.5,.4);
\end{tikzpicture}
\caption{The three pairs of intervals approaching $\eta$ at each selected level. The ratios between successive lengths are not represented to scale.}
\label{fig:targets}
\end{figure}
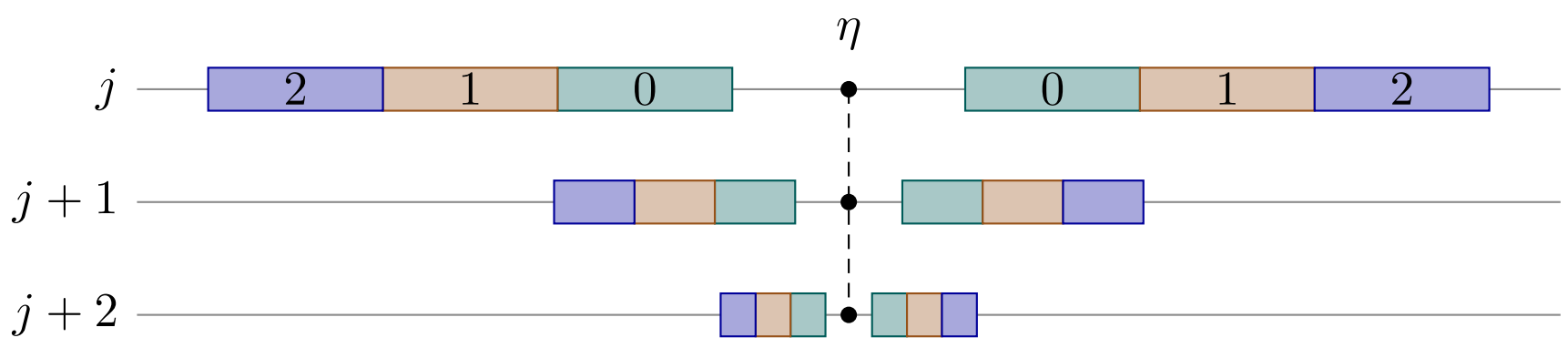

Put $U_i=\bigcup_jP_{j,i}$ and $U=\bigcup_{i=0}^2U_i$. Define $f_i$ to agree with $a$ on $U\setminus U_i$ and with the identity elsewhere; define $a'$ to agree with $a$ on $\mathcal{X}_\theta\setminus U$ and with the identity on $U$. These are involutions. They are continuous because the pieces are clopen away from $\eta$ and $a(\eta)=\eta$. We put
\[
 G_\theta=\langle a',f_0,f_1,f_2,b\rangle.
\]

\begin{proposition}\label{prop:periodic}
The group $G_\theta$ is infinite and periodic, with the same orbits as the original dihedral group action.
\end{proposition}
\begin{proof}
The three involutions $f_i$ generate a Klein four-group as in the case of the Grigorchuk group (see \cite{Gri80,Gri84,Vor12,Nek18}), and $a'$ commutes with them since it acts on a disjoint set. The finite group $A=\langle a',f_0,f_1,f_2\rangle$ is a fragmentation of $a$ in the sense of \cite{Nek18}. Every element of $A$ fixes one of the sets $U_i$ pointwise, and every $U_i$ accumulates at $\eta$. Periodicity follows from \cite[Theorem~1.1]{Nek18}, taking the other fragmentation to be $\langle b\rangle$. Every edge of the original orbital graph joining distinct vertices is retained by some generator. Thus, the orbits are unchanged and the group is infinite. The $f_i$ are distinct involutions since all three families of intervals are nonempty. The involution $a'$ is nontrivial on the gaps between successive families, and $b$ is nontrivial.
\end{proof}

\section{Realization by a tile inflation process}\label{sec:native}
Using the conventions of \subsecref{subsec:prelim-tiles} (see \cite[Section~2]{Kua25} and \cite{Kua26a}), we realize the fragmentation over the diagram of \cite[Subsection~5.1]{Kua25}. Write $\mathfrak H_\theta$ for the groupoid of germs of the fragmented action.

\subsection{The Rokhlin--Kakutani partitions}
Let $\theta=[0;c_1,c_2,\ldots]$, with convergents $p_m/q_m$, and put $\beta_m=|q_m\theta-p_m|$, $\beta_{-1}=1$. In the geometric realization of the Rokhlin--Kakutani partitions $\mathcal A_m$ of \cite[Theorem~5.2]{Kua25}, there are two towers at each level. Their bases, in the coordinates of the original circle, are the intervals of lengths $\beta_m$ and $\beta_{m-1}-\beta_m$ in $[0,\beta_{m-1}]$, with the one-sided conventions for the split endpoints. The incidence matrix between consecutive levels, with rows indexed by the towers at the higher level, is
\begin{equation}\label{eq:native-incidence}
 M(c)=\begin{pmatrix}1&c\\1&c-1\end{pmatrix}.
\end{equation}
The partitions are obtained by the following induction. In the normalized base at a step with angle $t\in(0,1)$ and $c=\lfloor1/t\rfloor$, take the intervals
\[
 I_i=[it,1-(c-i)t]\quad(0\leq i\leq c),\qquad
 J_i=[1-(c-i)t,(i+1)t]\quad(0\leq i<c).
\]
We use the appropriate one-sided endpoints for these clopen intervals. The maps $A(s)=1-s$ and $B(s)=t-s\pmod1$ give two paths, in the orders
\[
 I_0,I_c,I_1,I_{c-1},\ldots,
 \qquad J_0,J_{c-1},J_1,J_{c-2},\ldots,
\]
where each vertex is used once and the edges have alternating labels $A,B$. On the rescaled base $[0,1]$, the unused $B$-edges at the bases induce the map $s\mapsto1-s$, and the unused edges at the other endpoints, transported to the bases, induce the map $s\mapsto t'-s\pmod1$, where $t'=1/t-c$, with the corresponding one-sided conventions at the split points (see \cite[Proposition~5.1]{Kua25}). This follows by substituting the endpoints of $I_i,J_i$ in the two reflections: the unused generator maps the interval at the base into the union of the two bases. If $c=1$, then $I_0,I_1$ give a path of two vertices and $J_0$ consists of a single vertex, with both $A$ and $B$ giving boundary edges at its only vertex. The same induced maps apply in this case. Replacing the vertices of these paths by copies of the lower-level tiles and joining their boundary edges gives the tile inflation of \cite[Subsection~5.1]{Kua25}. The internal edges form a path with alternating labels.

The lengths of the bases at level $m+1$ are $\beta_{m-1}-c_{m+1}\beta_m=\beta_{m+1}$ and $\beta_m-\beta_{m+1}$. Each floor is obtained by a composition of the two original isometries, so its projected length is at most $\beta_{m-1}$. The heights $h_{0,m},h_{1,m}$ of the two towers have initial values $(1,1)$ and satisfy the recursion given by the incidence matrix in \eqref{eq:native-incidence}. It follows that
\begin{equation}\label{eq:towersize}
 h_{0,m}=q_m+q_{m-1},\qquad h_{1,m}=q_m,
 \qquad \beta_{m-1}<1/q_m\quad(m\geq1).
\end{equation}
The last inequality is the standard estimate for convergents. The partitions generate the topology: by Lemma~\ref{lem:endpoints}, every split endpoint eventually occurs as an endpoint of a floor, and the maximum projected length tends to zero. Write $\Omega_n(\mathsf B)$ for the finite paths through level $n$. By \subsecref{subsec:prelim-tiles}, nested floors identify $\Omega(\mathsf B)$ with $\mathcal{X}_\theta$, and the equivalence relations within the towers give the tail groupoid $\mathfrak T_{\mathsf B}$.

\begin{figure}[htbp]
\centering
\begin{tikzpicture}[x=.82cm,y=.85cm]
\foreach \x/\lab in {0/I_0,2/I_3,4/I_1,6/I_2}{\node[draw,rounded corners=1pt,minimum width=1.04cm,minimum height=.55cm] (I\x) at (\x,1.6) {$\lab$};}
\draw[conn,thick] (I0)--(I2) node[midway,above]{$A$};\draw[bound,thick] (I2)--(I4) node[midway,above]{$B$};\draw[conn,thick] (I4)--(I6) node[midway,above]{$A$};
\foreach \x/\lab in {0/J_0,2/J_2,4/J_1}{\node[draw,rounded corners=1pt,minimum width=1.04cm,minimum height=.55cm] (J\x) at (\x,0) {$\lab$};}
\draw[conn,thick] (J0)--(J2) node[midway,above]{$A$};\draw[bound,thick] (J2)--(J4) node[midway,above]{$B$};
\foreach \name/\lab in {I0/B,I6/B,J0/B,J4/A}{\draw[dashed,black!60] (\name.south)--++(0,-.38) node[below]{$\lab$};}
\node[align=left,right] at (7.1,.8){Dashed edges:\\unused boundary edges};
\end{tikzpicture}
\caption{The paths in the two towers at a normalized induction step with $c=3$. At later levels the vertices are replaced by copies of the lower-level tiles.}
\label{fig:paths}
\end{figure}
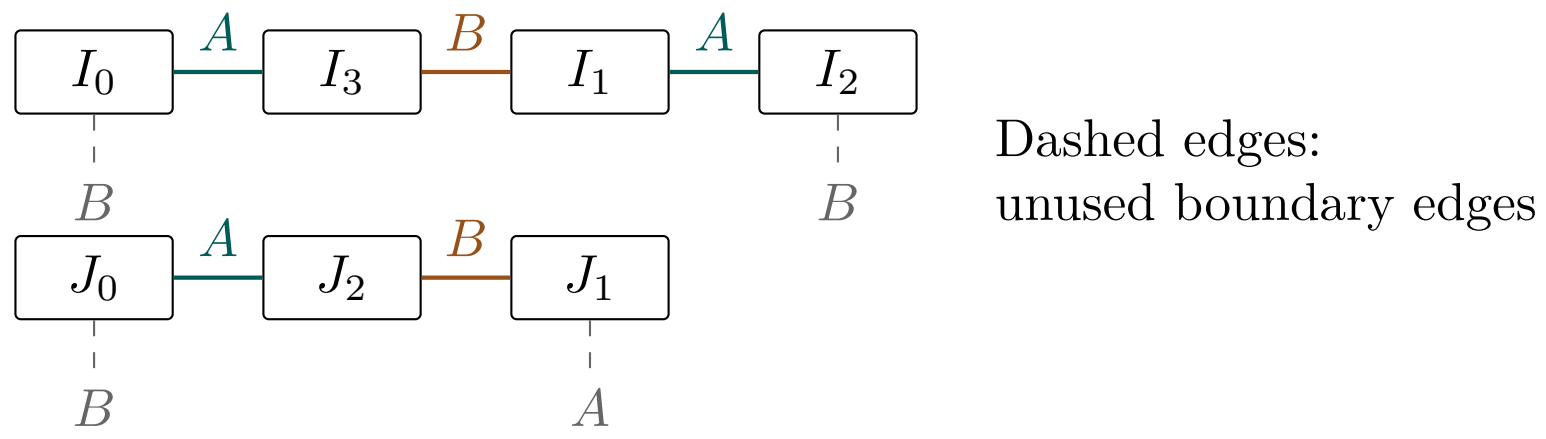

\subsection{Endpoints and boundary points}\label{subsec:nativeboundary}
Let $E_m^\circ$ denote the projected endpoints of the floors of $\mathcal A_m$.
\begin{lemma}\label{lem:endpoints}
For $m\geq1$,
\[
 \{k\theta\pmod1:|k|\leq\lfloor(q_m-3)/4\rfloor\}
 \subseteq E_m^\circ\subseteq
 \{k\theta\pmod1:|k|\leq3q_m\}.
\]
\end{lemma}
\begin{proof}
The endpoints of the bases at level $m$ are $0,\beta_m,\beta_{m-1}$ in the original coordinates, so they have the form $k\theta\pmod1$ with $|k|\leq q_m$. Each floor is obtained from its base by at most $2q_m-1$ steps along the corresponding finite path, by \eqref{eq:towersize}. On $k\theta$ the maps act by $k\mapsto-k$ and $k\mapsto1-k$, so the upper bound follows.

The base with endpoint $0$ lies in a path of at least $q_m$ vertices. Starting at this base, one direction along the path has length at least $(q_m-1)/2$. The alternating maps on the integers $k$, applied to $0$, produce successive terms of $0,1,-1,2,-2,\ldots$, with a possible initial repetition of $0$ or a reversal of signs. All these terms give projected endpoints of floors of the tower, proving the lower bound.
\end{proof}

For $s\in\{a,b\}$, let $U_m(s)$ be the union of floors corresponding to vertices whose $s$-edges are boundary edges at level $m$. The set $U_m(s)$ is $s$-invariant since the internal edges labeled by $s$ occur in inverse pairs.
\begin{lemma}\label{lem:boundary}
For $m$ such that the projected lengths of the floors are less than $1/2$,
\[
 U_m(s)\subseteq\{\zeta:\operatorname{dist}(\zeta,\operatorname{Fix}(s))\leq2\beta_{m-1}\},
\]
where $\operatorname{Fix}(s)\subseteq\mathcal X_\theta$ and distances are measured after projection to the circle.
\end{lemma}
\begin{proof}
There are four pairs $(v,s)$, where $v$ is a vertex of one of the two paths and $s\in\{a,b\}$ is the label of a boundary edge at $v$. A path consisting of a single vertex contributes two pairs. Every floor containing a fixed point of $s$ corresponds to such a pair, since an internal edge joins distinct floors. Since each floor has length less than $1/2$, the two fixed points of $b$ belong to distinct floors. Thus, three of the four pairs correspond to floors containing fixed points of the respective generators. There is at most one other pair. All the floors containing the fixed points satisfy the assertion, by \eqref{eq:towersize}.

Let $C$ be an additional floor in $U_m(s)$. If $s(C)$ does not meet the other floors of $U_m(s)$, then invariance gives $s(C)\subseteq C$, and hence $s(C)=C$, since $s$ is an involution. If a proper interval of length less than $1/2$ is invariant under a reflection, then it contains a fixed point of the reflection in its interior. The fixed point $0$ of the reflection on the circle cannot lie in that interior, since $0$ is an endpoint of the partition. Thus, $C$ would contain an actual fixed point in the split circle, a contradiction. It follows that $s(C)$ meets a floor containing a fixed point of $s$. Both floors have lengths at most $\beta_{m-1}$, and reflection preserves the distance to its fixed points. This proves the assertion.
\end{proof}

\subsection{Labeling the connectors}
\begin{theorem}\label{thm:hierarchy}
For every irrational $\theta\in(0,1)$, the action of $G_\theta$ on the split circle $\mathcal{X}_{\theta}$ is conjugate to an action defined by a tile inflation process over a telescoping of the Bratteli diagram $\mathsf B$ constructed in \cite[Theorem~5.2]{Kua25} for the associated dihedral group action. The tail groupoid $\mathfrak T_{\mathsf B}$ in the realization of the dihedral group action is an open subgroupoid of the groupoid of germs of $G_\theta$, with discrete complement.
\end{theorem}

Put $m_j=n_j+5$, where $n_j$ is as in \secref{sec:construction}. We show that the $a$-edges with endpoints in $P_{j,i}$ are labeled in the inflation from level $m_{j-1}$ to level $m_j$ of the telescoped diagram.

Since $q_{m_j}\geq8Q_j$, Lemma~\ref{lem:endpoints} shows that all projected endpoints of $P_{j,i}$ belong to $E_{m_j}^\circ$. By the one-sided conventions, $P_{j,i}$ is a union of cylinders of length $m_j$. By Lemma~\ref{lem:boundary}, every point of $U_{m_j}(a)$ is within distance $2\beta_{m_j-1}<1/(4Q_j)$ of $\eta$, whereas the points of $P_{j,i}$ have distance at least $10/Q_j$ from $\eta$. Thus, the corresponding edges labeled by $a$ are internal edges of the tiles at level $m_j$.

To see that the labels of internal edges are preserved, recall that $D_j=q_{n_j+12}\geq21q_{m_j}$ and $Q_{j+1}\geq10946D_j$. All projected endpoints of the floors at level $m_j$ have the form $k\theta\pmod1$ with $|k|\leq3q_{m_j}$. For $k\ne0$, applying best approximation to $2k$ gives
\[
 \operatorname{dist}(\eta,k\theta)>1/(4D_j),
\]
since $2|k|<D_j$ and $\|2k\theta\|>1/(2D_j)$. The case $k=0$ is immediate. Every $P_{k,i}$ with $k>j$ lies within distance $50/Q_{j+1}<1/(200D_j)$ of $\eta$. Consequently, both endpoints of an $a$-edge to be labeled at a later level lie in the floor of level $m_j$ containing $\eta$. An internal edge joins distinct floors, so this edge cannot have been introduced by level $m_j$.

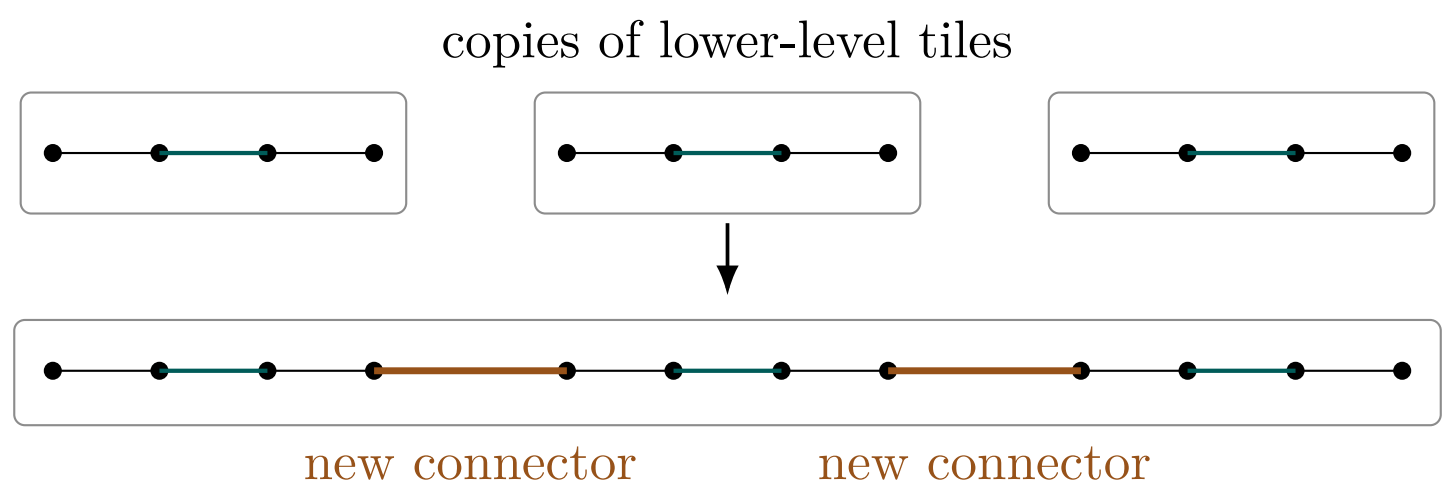
\begin{figure}[htbp]
\centering
\begin{tikzpicture}[x=.85cm,y=.8cm]
\foreach \x in {0,4,8}{
 \draw[black!45,rounded corners=2pt] (\x-.25,1.1) rectangle (\x+2.75,2.1);
 \draw (\x,1.6)--(\x+2.5,1.6);
 \foreach \d in {0,.83,1.67,2.5}{\node[point] at (\x+\d,1.6) {};}
 \draw[conn,thick] (\x+.83,1.6)--(\x+1.67,1.6);
}
\node at (5.25,2.5){copies of lower-level tiles};
\draw[edge] (5.25,1.02)--(5.25,.42);
\draw[black!45,rounded corners=2pt] (-.3,-.65) rectangle (10.8,.22);
\foreach \x in {0,4,8}{\draw (\x,-.2)--(\x+2.5,-.2);\foreach \d in {0,.83,1.67,2.5}{\node[point] at (\x+\d,-.2) {};}\draw[conn,thick] (\x+.83,-.2)--(\x+1.67,-.2);}
\draw[bound,line width=1.3pt] (2.5,-.2)--(4,-.2) (6.5,-.2)--(8,-.2);
\node[below,bound] at (3.25,-.65){new connector};\node[below,bound] at (7.25,-.65){new connector};
\end{tikzpicture}
\caption{The labels of the internal edges are preserved when the labels of the newly introduced connectors are assigned.}
\label{fig:connectors}
\end{figure}

Replace each pair of boundary edges labeled by $a$ by a pair for each of the four fragments of $a$, and retain the boundary edges labeled by $b$. At level $m_j$, a newly introduced $a$-edge whose corresponding floor is contained in $P_{j,i}$ receives the labels $f_\ell$ for $\ell\ne i$, and every other newly introduced $a$-edge receives the label $a'$. The remaining labels of fragments of $a$ give loops at the two endpoints by joining the corresponding boundary edges at each endpoint. Use the same labeling of internal edges, loops and boundary edges at the initial level, retaining all labels $b$. This is the convention for graphs with boundary in \cite[Subsection~2.1]{Kua25}. The labels assigned at earlier levels are preserved. By Lemma~\ref{lem:boundary}, every original edge joining distinct vertices eventually appears in a finite tile, and thus the infinite tiles have exactly the edges between distinct vertices of the fragmented action. At a fixed point of an original reflection, all its fragments act trivially. The homeomorphisms obtained by extending these maps continuously are those of \secref{sec:construction}, proving the conjugacy in Theorem~\ref{thm:hierarchy}.

The number of floors contained in $P_{j,i}$ is comparable to $q_{m_j}/Q_j$. The lower bound follows from the lengths of the floors in \eqref{eq:towersize} and of the intervals defining $P_{j,i}$; the upper bound follows from Lemmas~\ref{lem:endpoints} and~\ref{lem:arithmetic}, using the lower bound on the distances between $q_{m_j}$ consecutive iterates. Thus, the number of connectors labeled by fragments of $a$ need not be uniformly bounded.

\subsection{The tail groupoid and the groupoid of germs}\label{subsec:commontail}
For the realization of the original dihedral group action, the tail groupoid is generated by the germs $(a,\zeta)$ with $a(\zeta)\ne\zeta$ and $(b,\zeta)$ with $b(\zeta)\ne\zeta$. Indeed, Lemma~\ref{lem:boundary} shows that each of these germs is eventually represented by an internal edge of a finite tile. Conversely, each replacement of a finite prefix is realized by the unique path of internal edges between two floors of a tower. Thus, the tail groupoid is principal and has the same orbits as the original dihedral group action. By construction, each of these germs is realized by a fragment of $a$ or by $b$. Hence, $\mathfrak T_{\mathsf B}\subseteq\mathfrak H_\theta$ and the tail groupoid is the same subgroupoid for both actions.

Let $E=\operatorname{Fix}(a)\cup\operatorname{Fix}(b)$. If $s$ is a generator of $G_\theta$ and $\zeta\notin E$, then $(s,\zeta)\in\mathfrak T_{\mathsf B}$ since it is one of the germs just described or the identity germ. For a word $s_m\cdots s_1$, its bisection can contain an arrow outside the tail groupoid only at a point of the finite set
\[
 \bigcup_{i=0}^{m-1}(s_i\cdots s_1)^{-1}(E),
\]
where the word for $i=0$ is the identity. These bisections form a basis for the topology of $\mathfrak H_\theta$, and the source map is injective on each of them. The complement of the tail groupoid is therefore discrete. This completes the proof of Theorem~\ref{thm:hierarchy}.

For comparison with the rotation action, consider the orientation-preserving part of the tail groupoid for the original dihedral group action. Let $\mathfrak T^{\rm rot}_\theta=\mathfrak T_{\mathsf B}\cap\mathfrak N_\theta$, where $\mathfrak N_\theta$ is identified with the rotation subgroupoid of $\mathfrak D_\theta$. It is an AF subgroupoid, obtained by intersecting the elementary subgroupoids with the kernel of the orientation cocycle. The rotation by $+\theta$ is $r_\theta=ba$. Exactly one step of $ba$ fixes its starting point if and only if the initial point belongs to
\[
 E'=\{1/2,-\theta/2,(1-\theta)/2\}\pmod1.
\]
At such a point, the rotation arrow differs from the unique arrow of $\mathfrak T_{\mathsf B}$ with the same source and range by a nontrivial isotropy arrow. Otherwise, both steps move their starting points. The points of $E'$ lie in distinct rotation orbits. A path of $r_\theta^m$ uses at most one of the three arrows $(r_\theta,\zeta)$, $\zeta\in E'$, or their inverses. Its germ belongs to $\mathfrak T_{\mathsf B}$ precisely when it avoids these arrows. Thus, $\mathfrak T^{\rm rot}_\theta$ is obtained by deleting the three arrows from the generating bisection of $r_\theta$. In particular, it is not the AF subgroupoid obtained by the backward modular restrictions in \subsecref{subsec:compression}. The latter is obtained by deleting the arrows of $S_\theta$ entering a set of one or two points. The tail groupoid of the Bratteli diagram for the dihedral group action is not contained in the rotation groupoid.

At $\eta$, the germ of the original reflection $a$ does not belong to $\mathfrak H_\theta$. Indeed, the orbital graph of $\eta$, with loops deleted, is a ray. In a closed walk based at $\eta$, the loops at other vertices give trivial germs. The edges between distinct vertices cancel in inverse pairs along the ray. Thus, its germ belongs to the group of germs of $A$ at $\eta$. The three sets $U_i$ and the gaps between the intervals at successive levels accumulate at $\eta$, so the group of germs at $\eta$ in $\mathfrak H_\theta$ is isomorphic to $A\cong(\mathbb Z/2\mathbb Z)^3$. Every element of $A$ fixes one of the $U_i$ pointwise, whereas the original reflection moves all their points. Consequently, no such germ agrees with $(a,\eta)$. The groupoids of germs of the two actions contain the same tail groupoid. However, the groupoid of germs of the dihedral group action is not contained in that of the fragmented action. At the other two fixed points of the reflections, the groups of germs still have order two, since all trajectories in their orbits avoid $\eta$, where the fragments do not all agree locally with elements of the dihedral group.

\subsection{Unbounded incidence matrices}
By \cite[Corollary~5.6]{Kua25}, the original action is of bounded type if and only if the partial quotients are bounded. Each incidence matrix $M(c_n)$ in \eqref{eq:native-incidence} has no zero rows or columns and contains the entry $c_n$. If we multiply on either side by nonnegative integer matrices with no zero rows or columns, then the resulting matrix has an entry at least $c_n$. If the $c_n$ are unbounded, then every telescoping has unbounded incidence matrices. On the other hand, the two finite tiles have at most four boundary vertices in total at each level. By Lemma~\ref{lem:boundary}, the boundary points of the infinite tiles are precisely the three fixed points of the reflections. The bound on the numbers of edges of the diagram in \cite[Definition~2.13]{Kua25} is not satisfied by any such telescoping and is not used in the growth estimate.

\section{Contraction of traverses}\label{sec:contraction}
Choose $x$ with trivial stabilizer in the dihedral group action on $\mathcal X_{\theta}$ whose projection back to the circle is not in $\mathbb Z\theta\pmod1$. Its linear orbital graph has consecutive vertices $v_{2k}=x+k\theta$ and $v_{2k+1}=-x-k\theta$. The fragmentation retains every edge between distinct vertices, so the distances are those in the line. The labels are determined by the sequence $\kappa(x+k\theta)$, $k\in\mathbb Z$, where $\kappa(x)\in\{0,1,2,*\}$ is the color of the $a$-edge at $x$ and the value $*$ corresponds to points outside $U$. We have $\kappa(-x)=\kappa(x)$ away from the split points and $\eta$. The orbit is dense, so the action on it is faithful. Each finite word of colors in this graph is realized on an open interval of projected points, since its finitely many arguments avoid $\eta$ and the split points. For a split point with trivial stabilizer in the dihedral group action, finite words are the corresponding one-sided limits, since the rotation preserves the sides of split points. Conversely, every word realized on an open interval occurs in each such orbit by density. Hence, the sequences of colors of all these linear orbital graphs have the same finite subwords, and we may count them on the chosen graph.

The orbital graph of $\eta$ is a ray whose consecutive vertices, starting at its endpoint, are
\[
 \eta,\ \eta+\theta,\ \eta-\theta,\ \eta+2\theta,\ \eta-2\theta,\ldots.
\]
Let $I_j$ be the segment of the first $16Q_j$ vertices and $p_j$ the next vertex. Put $V=(\C)^2$, with the three nonzero elements $u_i$ corresponding to $f_i$, $i=0,1,2$, and let $\chi_i:V\to\C$ have kernel $\{0,u_i\}$. For a path in this graph, we take the sum of the elements $u_i\in V$ corresponding to applications of $f_i$ at $\eta$; all other steps contribute zero. We count the return paths according to this sum rather than the germs of the corresponding words.

\begin{lemma}\label{lem:fold}
Consider an $a$-edge with endpoints in $P_{j+1,i}$ and the two branches of length $16Q_j$ extending away from it. Identifying the two branches, with the endpoints of the $a$-edge sent to $\eta$, gives a map onto $I_j\cup\{p_j\}$ preserving the labels of steps within the segment. The image of a complete crossing from one outer endpoint to the other has sum $\sigma$ in $V$ satisfying $\chi_i(\sigma)=1$.
\end{lemma}
\begin{proof}
The central edge has endpoints $\eta\pm t$ with $|t|\leq50/Q_{j+1}<1/(200D_j)$. Replacing $t$ by zero maps both branches onto $I_j\cup\{p_j\}$. At the edges labeled by $a$ other than the central edge, the values of $\kappa$ agree: their coordinates are $\eta+k\theta\pm t$, with nonzero $k$ of absolute value at most $16Q_j+1$. The interval endpoints through level $j$ project to $\ell\theta\pmod1$ with $|\ell|\leq Q_j$. Since $|k-\ell|\leq17Q_j+1<100Q_j$, equation~\eqref{eq:sep2} ensures that none of these endpoints lies between the coordinates $\eta+k\theta$ and $\eta+k\theta\pm t$. By \eqref{eq:sep1}, the distance from $\eta+k\theta\pm t$ to $\eta$ is greater than $1/(2D_j)-|t|>50/Q_{j+1}$. Thus, no vertex other than the endpoints of the central edge enters the neighborhood containing the pieces at later levels. The edges labeled by $b$ are unchanged. The two generators $f_\ell$ for which $\chi_i(u_\ell)=1$ interchange the central endpoints; the other fragments of $a$ fix them. In the orbital graph at $\eta$, all fragments of $a$ fix $\eta$. A path from one outer endpoint to the other crosses the central edge an odd number of times, and hence $\chi_i(\sigma)=1$. The same argument applies to such a segment in the orbital graph of a fixed point of a reflection, provided that it does not contain the endpoint of the ray.
\end{proof}

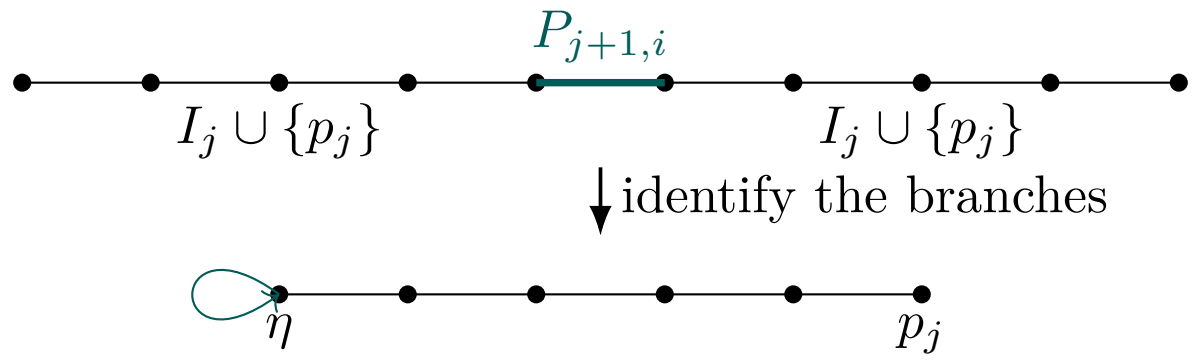
\begin{figure}[htbp]
\centering
\begin{tikzpicture}[x=.85cm,y=.7cm]
\draw (0,2)--(9,2); \foreach \x in {0,1,2,3,4,5,6,7,8,9} \node[point] at (\x,2) {};
\draw[conn,line width=1.4pt] (4,2)--(5,2); \node[above,conn] at (4.5,2) {$P_{j+1,i}$};
\node[below] at (2,2) {$I_j\cup\{p_j\}$}; \node[below] at (7,2) {$I_j\cup\{p_j\}$};
\draw[edge] (4.5,1.2)--(4.5,.55); \node[right] at (4.5,.9) {identify the branches};
\draw (2,0)--(7,0); \foreach \x in {2,3,4,5,6,7} \node[point] at (\x,0) {};
\draw[conn,->] (2,0) .. controls (1.1,.8) and (1.1,-.8) .. (2,0);
\node[below] at (2,0) {$\eta$}; \node[below] at (7,0) {$p_j$};
\end{tikzpicture}
\caption{The map identifying the two branches at an $a$-edge with endpoints in $P_{j+1,i}$ with the segment $I_j\cup\{p_j\}$ of the orbital graph at $\eta$.}
\label{fig:fold}
\end{figure}

Let $w=s_m\cdots s_1$ be a word in the five involutions. For each index $k$, we start at $p_j$ immediately before applying $s_k$. Let $T_{j,\sigma}(w)$ count those indices $k$ for which the path enters $I_j$ on this step, returns to $p_j$ before or at the end of the word, and has sum $\sigma\in V$ up to and including its first return to $p_j$. The initial time determines the return path, which ends at the first return to $p_j$. Put $T_j(w)=\sum_{\sigma\ne0}T_{j,\sigma}(w)$.

\begin{lemma}\label{lem:contraction}
For every word $w$, $T_{j+1}(w)\leq(2/3)T_j(w)$. In particular,
\begin{equation}\label{eq:contraction}
 T_j(w)\leq |w|(2/3)^{j-1}.
\end{equation}
\end{lemma}
\begin{proof}
By Lemma~\ref{lem:arithmetic}, among $2Q_{j+1}\leq k<3Q_{j+1}$ there is an edge with endpoints in $P_{j+1,i}$ for each $i$. For each $i$, fix one such edge. In the path metric of the ray, these edges lie between distances $4Q_{j+1}-1$ and $6Q_{j+1}$ from $\eta$, so the segments described in Lemma~\ref{lem:fold} are contained in the interior of $I_{j+1}$. Every return path counted by $T_{j+1}$ reaches $\eta$ and crosses the three segments.

Fix $i$. Before the return path at level $j+1$ first reaches $\eta$, take the traverse towards $\eta$ of the corresponding segment, as in \subsecref{subsec:prelim-traverses}. Since both endpoints are identified with $p_j$ by Lemma~\ref{lem:fold}, its image is a first-return path counted by $T_{j,\sigma}$ with $\chi_i(\sigma)=1$. The starting time of this traverse and the fixed entrance vertex determine its position in the larger return path. Reversing the letters of $w$ from this vertex and time until the first visit to $p_{j+1}$ recovers the initial time of the return path at level $j+1$. Its final time is determined by the first return to $p_{j+1}$. Thus, $T_{j+1}(w)\leq\sum_{\chi_i(\sigma)=1}T_{j,\sigma}(w)$. Summing over $i=0,1,2$, we count each nonzero $\sigma$ twice. This gives $3T_{j+1}\leq2T_j$ by specializing the counting argument of \cite[proof of Theorem~3.1]{BNZ22} to the three nonzero characters of $V$. Finally, $T_1(w)\leq |w|$.
\end{proof}

\section{Counting finite configurations and portraits}\label{sec:portraits}
\subsection{Counting labeled balls and neighborhoods of edges}
Let $p(r)$ count the isomorphism classes of oriented rooted labeled balls of radius $r$ in the chosen linear orbital graph. Let $A_j(R)$ count the isomorphism classes of neighborhoods of radius $2R+2$ of $a$-edges with endpoints in some $P_{j,i}$, with both endpoints distinguished and the orientation retained.

\begin{lemma}\label{lem:complexity}
There is a constant $C$ such that
\begin{equation}\label{eq:complexity}
 p(r)\leq C(r+1),\qquad A_j(R)\leq C(1+R/Q_j)\quad (R\geq Q_j).
\end{equation}
\end{lemma}
\begin{proof}
Consider the sequence $\kappa(x+k\theta)$ for $|k|\leq M$, with $M=r+2$ for the rooted balls and $M=2R+2$ for the neighborhoods of $a$-edges. Choose $m$ with $Q_m\leq M<Q_{m+1}$ (if $M<Q_1$, then we omit the first collection of endpoints below and use $1$ in place of $m+1$). The projected endpoints of $P_{j,i}$ for $j\leq m$, translated by $-k\theta$, have the form $\ell\theta\pmod1$ with $|\ell|=O(M)$. Each of the at most eight projected endpoints at level $m+1$ contributes $2M+1$ translates. Together, these give $O(M+1)$ points.

For subsequent levels, the translated endpoints lie in neighborhoods of $\eta-k\theta$, $|k|\leq M$, of radius $50/Q_{m+2}$. Equations~\eqref{eq:sep1}--\eqref{eq:sep2} with $j=m+1$ and \eqref{eq:fib} imply that these neighborhoods are pairwise disjoint and contain none of the points counted above. In each such neighborhood, at most one entry of the sequence of colors varies, giving at most four words. This proves the first bound, since the root has only a bounded number of positions relative to its $a$-edge.

For the second bound, restrict $x$ to the projections of the two intervals forming $P_{j,i}$. Summing over $i$ changes only the constant. Each interval has length $O(1/Q_j)$. By Lemma~\ref{lem:arithmetic}, each block of $Q_j$ consecutive iterates meets such an interval in a bounded number of points. Dividing the segments of the orbit of length $O(R)$ into such blocks bounds the number of points from levels through $m+1$ in the interval by $O(1+R/Q_j)$. Apply the same argument to the points $\eta-k\theta$, enlarging the interval by $50/Q_{m+2}$ at each end. The enlargement is still of length $O(1/Q_j)$. Each neighborhood gives at most four words, as above. This proves the second bound.
\end{proof}

\subsection{The portrait estimate}
We use the standard portraits recalled in \subsecref{subsec:prelim-traverses} on the linear orbital graph chosen in \secref{sec:contraction}, on which the action is faithful.

\begin{lemma}\label{lem:portrait}
For $j\geq2$ and $R\geq16Q_j$, every word $w$ of length at most $R$ satisfies
\begin{equation}\label{eq:portrait}
 |\mathcal P_w|\leq C\left(Q_j+\frac{R^2}{Q_j}(2/3)^j\right).
\end{equation}
\end{lemma}
\begin{proof}
For trajectories of radius less than $16Q_j$, the corresponding elements of the portrait are determined by the oriented rooted balls of radius $16Q_j$ and the word $w$. By Lemma~\ref{lem:complexity}, there are at most $CQ_j$ such elements.

A trajectory of radius at least $16Q_j$ covers the segment from its initial vertex to a vertex visited at the maximum distance from it. Remove segments of length $16Q_{j-1}+2$ from both ends. By \eqref{eq:fib}, the remaining segment has length greater than $4Q_j$ and hence contains a block of $Q_j$ consecutive $a$-edges. The projection of the right-hand interval forming $P_{j,0}$ has length at least $8/Q_j$. Lemma~\ref{lem:arithmetic} applied to the coordinates of this block of edges gives an $a$-edge with endpoints in $P_{j,0}$. Since this edge lies in the remaining segment, its two branches of length $16Q_{j-1}$ in Lemma~\ref{lem:fold} lie in the original segment. Orient the segment consisting of this edge and its branches in the direction from the initial vertex to the chosen vertex. Take a traverse as in \subsecref{subsec:prelim-traverses}. The two endpoints are identified with $p_{j-1}$, and no interior vertex maps to $p_{j-1}$. Thus, by Lemma~\ref{lem:fold}, its image is a return path counted by $T_{j-1}(w)$.

Record the oriented labeled neighborhood of radius $2R+2$ of the central edge, together with the direction and starting time of the traverse. The direction specifies the entrance vertex, from which the preceding letters of the word, inverted and applied in reverse order, recover the root. Since the trajectory visits the central edge, both endpoints of this edge are at distance at most $R+1$ from the root. Thus, the recorded neighborhood contains the ball representing the corresponding element of the portrait. The recorded data determine this element and have at most $2A_j(R)T_{j-1}(w)$ possible values. Using \eqref{eq:contraction} and \eqref{eq:complexity} gives the second term in \eqref{eq:portrait}.
\end{proof}

\Needspace{20\baselineskip}
\section{Growth of the groups and repetitivity of their orbital graphs}\label{sec:repetitivity}
\begin{theorem}\label{thm:growth}
For each irrational $\theta\in(0,1)$, the fragmentation of the associated minimal dihedral group action constructed in \secref{sec:construction} gives an infinite periodic group $G_\theta$ generated by $5$ involutions. Let $\gamma_{G_\theta}$ be its growth function with respect to these generators, and let $q_n(\theta)$ be the denominators of the continued-fraction convergents of $\theta$. If
\[
 \liminf_{n\to\infty}(2/3)^{n/64}\log q_n(\theta)=0,
\]
then $G_\theta$ has intermediate growth. Put
\[
 \Lambda_\theta=\limsup_{n\to\infty}\frac{\log q_n(\theta)}{n}.
\]
If $\Lambda_\theta<\infty$, then for every
\[
 \frac{64\Lambda_\theta}{64\Lambda_\theta+\log(3/2)}<\alpha<1,
\]
there exists $C>0$ such that $\gamma_{G_\theta}(R)\leq\exp(CR^\alpha)$. For every irrational $\theta$ and every point with trivial stabilizer in the dihedral group action, the number of isomorphism classes of rooted labeled balls of radius $r$ in the orbital graph of the fragmented action at this point is $O(r+1)$. If the partial quotients of $\theta$ are unbounded, then the orbital graphs at points with trivial stabilizer in the dihedral group action are not linearly repetitive.
\end{theorem}

By Proposition~\ref{prop:periodic}, the group is infinite and periodic. Lemma~\ref{lem:complexity} gives the bound on the number of isomorphism classes of rooted labeled balls. The remaining assertions are proved in the next two subsections.

\subsection{The upper bound on growth}\label{subsec:uppergrowth}
Put $\delta=\log(3/2)$ and $R_j=\lceil Q_j\exp(\delta j/2)\rceil$. For all sufficiently large $j$, Lemma~\ref{lem:portrait} gives $|\mathcal P_w|\leq CQ_j$ for all words of length at most $R_j$. The number of isomorphism classes of bi-rooted oriented labeled balls of radius at most $R_j$ is at most
\[
 \sum_{r=0}^{R_j}p(r)(2r+1)=O(R_j^3).
\]
Write $N_j$ for this number and $M_j=\lceil CQ_j\rceil$ for an upper bound on the cardinality of a portrait. There are at most $(N_j+1)^{M_j}$ possible portraits. Indeed, we can list the balls of each portrait in a fixed order and fill the remaining entries with an extra symbol. Since equal portraits determine equal elements of $G_\theta$ by \subsecref{subsec:prelim-traverses}, the counting argument of \cite[Proposition~2.3]{BNZ22} yields
\begin{equation}\label{eq:growth}
 \log\gamma_{G_\theta}(R_j)\leq CQ_j\log R_j.
\end{equation}
Dividing \eqref{eq:growth} by $R_j$ gives
\[
 \frac{\log\gamma_{G_\theta}(R_j)}{R_j}
 \leq Ce^{-\delta j/2}(\log Q_j+j).
\]
If $\liminf_n e^{-\delta n/64}\log q_n=0$, then $\liminf_j e^{-\delta j/2}\log Q_j=0$. Indeed, for $n_j\leq n<n_{j+1}$, we have $q_{n_j}\leq q_n$ and $n=32j+O(1)$. Thus, the quotient $\log\gamma_{G_\theta}(R_j)/R_j$ tends to zero along a subsequence. By submultiplicativity, $\lim_{R\to\infty}R^{-1}\log\gamma_{G_\theta}(R)=\inf_{R\geq1}R^{-1}\log\gamma_{G_\theta}(R)$, where the radii are integers. Hence, the group has subexponential growth. Since it is infinite and periodic, it has intermediate growth by Gromov's theorem \cite{Gro81} (see \subsecref{subsec:prelim-growth}).

If $\Lambda_\theta<\infty$, then $\log R_j=O(j)$ and $\limsup_j j^{-1}\log R_j\leq32\Lambda_\theta+\delta/2$. Fix $K>32\Lambda_\theta+\delta/2$. If $R$ is sufficiently large and $j=\lfloor (\log R)/K\rfloor$, then $R_j\leq R$. Submultiplicativity of the growth function and \eqref{eq:growth} imply
\[
 \log\gamma_{G_\theta}(R)
 \leq\lceil R/R_j\rceil\log\gamma_{G_\theta}(R_j)
 \leq C R j e^{-\delta j/2}.
\]
This is $O(R^\alpha)$ for every $\alpha>1-\delta/(2K)$. Letting $K$ tend to $32\Lambda_\theta+\delta/2$ gives the range of exponents in Theorem~\ref{thm:growth}.

\subsection{An obstruction to linear repetitivity}\label{subsec:nonlinear}
For the original action, the relation between linear repetitivity and bounded partial quotients was proved in \cite[Proposition~5.5]{Kua25}. Let $\mathcal R(\ell)$ be the least integer $L$ such that every block of length $L$ contains all words of length $\ell$ in the sequence of colors, up to reversal. By \subsecref{subsec:prelim-growth}, linear repetitivity implies $\mathcal R(\ell)=O(\ell)$, since consecutive $a$-edges have distance two in the line.

Let $p\ne\eta$ be an isolated discontinuity of $\kappa$ at which its two one-sided limits differ. There is an integer $k_0$ such that $\kappa(p+k_0\theta)\ne\kappa(\eta+k_0\theta)$ and both projected points are continuity points of $\kappa$. Indeed, take an interval on one side of $p$ on which $\kappa$ is constant. On the corresponding side of $\eta$ all three values $0,1,2$ occur arbitrarily close to $\eta$, so there is an open subinterval $J$ on which $\kappa(y)\ne\kappa(y+p-\eta)$, with both functions locally constant. By density, choose $\eta+k_0\theta\in J$. Put $d=\kappa(p+k_0\theta)$.

For large $n$, put $\epsilon_n=\norm{q_n\theta}$ and consider the set
\[
 E_n=\{x:\kappa(x)\ne\kappa(x+q_n\theta),\ \kappa(x+k_0\theta)=d\}.
\]
The set $E_n$ contains a one-sided interval of length $\epsilon_n$ at $p$, so its measure is positive. By the second condition, $E_n$ is disjoint from a fixed neighborhood of $\eta$. For large $n$, the short arc from $x$ to $x+q_n\theta$ stays outside a smaller neighborhood of $\eta$, where $\kappa$ has only finitely many discontinuities. If the values at the endpoints differ, then the arc contains one of these discontinuities. Hence, $m(E_n)\leq C_\theta\epsilon_n$.

Choose a point in the interior of the one-sided interval in $E_n$ whose iterates at indices from $\min(0,k_0)$ to $\max(q_n,k_0)$ are continuity points of $\kappa$. The corresponding word has length $\ell_n=q_n+O(1)$, and the initial point at index zero in each occurrence satisfies the conditions defining $E_n$. Since $\kappa$ is locally constant at every point of $\mathcal{X}_\theta$ other than $\eta$, the boundary of the set of points realizing this word is contained in the orbit of $\eta$ and has measure zero for $\mu_\theta$. By \subsecref{subsec:prelim-frag}, the frequency of this word is at most $C_\theta\epsilon_n$, even on orbits of points with trivial stabilizer in the dihedral group action that contain split points. The reversed word has the same frequency because the reflection preserves the colors and Lebesgue measure is invariant under reflections and rotations. Hence, $\mathcal R(\ell_n)\geq c_\theta/\epsilon_n$. Since $\epsilon_n<1/q_{n+1}$,
\begin{equation}\label{eq:rep}
 \frac{\mathcal R(\ell_n)}{\ell_n}\geq c_\theta\frac{q_{n+1}}{q_n}\geq c_\theta c_{n+1}.
\end{equation}
If the partial quotients are unbounded, then the orbital graph is not linearly repetitive. This completes the proof of Theorem~\ref{thm:growth}.

\subsection{Examples with unbounded partial quotients}
Let $\theta=[0;c_1,c_2,\ldots]$, where $c_1=3$, $c_n=2n$ for $n=2^m$, $m\geq1$, and $c_n=2$ otherwise. The partial quotients are unbounded and all but the first are even, so the angle does not satisfy the hypotheses of \cite[Subsection~5.3]{Kua25}.

Let $A(c)=\begin{pmatrix}c&1\\1&0\end{pmatrix}$. Since $A(2)\leq A(c)\leq (c/2)A(2)$ entrywise for $c\geq2$, the products satisfy
\[
 A(2)^n\leq A(c_1)\cdots A(c_n)
 \leq\left(\prod_{k=1}^n\frac{c_k}{2}\right)A(2)^n.
\]
We have $\sum_{k=1}^n\log(c_k/2)=O((\log n)^2)$, and the largest eigenvalue of $A(2)$ is $1+\sqrt2$. Consequently, the denominators of the convergents satisfy
\[
 \log q_n=n\log(1+\sqrt2)+O((\log n)^2).
\]
Theorem~\ref{main:growth} applies with $\Lambda_\theta=\log(1+\sqrt2)$. By \eqref{eq:rep}, the orbital graphs at points with trivial stabilizer in the dihedral group action are not linearly repetitive. Theorem~\ref{main:hierarchy} also applies, although no telescoping of the corresponding diagram has uniformly bounded numbers of edges at the levels.

For an example with $\Lambda_\theta=\infty$, take $c_n=\lceil e^n\rceil$ for all $n\geq1$. The inequalities $c_nq_{n-1}\leq q_n\leq(c_n+1)q_{n-1}$ give $\log q_n=n(n+1)/2+O(1)$. Thus, the first condition of Theorem~\ref{main:growth} holds, but the condition $\Lambda_\theta<\infty$ does not. The group has intermediate growth, and its orbital graphs at points with trivial stabilizer in the dihedral group action are not linearly repetitive.

\subsection{Open question} It remains to determine whether the groups in \secref{sec:construction} have intermediate growth for every irrational angle $\theta$.

\paragraph{AI disclosure.} The author used ChatGPT (OpenAI) and Claude (Anthropic) to assist with mathematical exploration, the development and examination of proofs, TikZ figures, and revision of the manuscript. The author is responsible for the mathematical content and the final text.

\end{document}